\documentclass[10pt,a4paper]{article}
\usepackage{makeidx}
\usepackage{amssymb}
\usepackage{graphics,graphicx}
\usepackage[T1]{fontenc}
\usepackage[francais]{babel}
\usepackage{alltt}

\newtheorem{hypo1}{Hypothèse}[section]

\newtheorem{thm1}[hypo1]{Théorème}

\newtheorem{hypo}{Hypothèse}[subsection]

\newtheorem{prop}[hypo]{Proposition}

\newtheorem{lem}[hypo]{Lemme}

\newtheorem{nota}[hypo]{Notation}

\begin{document}

\title{Vitesse dans le théorème limite central pour certains processus stationnaires fortement décorrélés}
 \date{\ }
        \author{Stéphane Le Borgne et Françoise Pène}

\maketitle

{\it Résumé : Nous prouvons un théorème limite central avec vitesse en $n^{-1/2}$ pour les processus stationnaires vérifiant une hypothèse de décorrélation forte. La démonstration  est une modification de la preuve du théorème de Rio. Elle est élémentaire mais longue et calculatoire. Nous montrons dans \cite{SLBF2} comment obtenir la décorrélation forte apparaissant dans l'hypothèse du théorème dans le cas de certains systèmes dynamiques quasi-hyperboliques.}

{\it Abstract : We prove a central limit theorem with speed  $n^{-1/2}$ for stationary processes satisfying a strong decorrelation hypothesis. The proof is a modification of the proof of a theorem of Rio. It is elementary but quite long  and technical.}

\section{Le théorème et sa démonstration}

Nous considérons un système dynamique probabilisé 
$\left(\Omega,{\cal F},\nu,T\right)$.
Nous notons ${\bf E}\left[\cdot\right]$ l'espérance relativement
à la mesure $\nu$~:
$${\bf E}[f]:=\int_\Omega f\, d\nu .$$
Pour toute fonction $f,g$ dans $L^2(\Omega,\nu)$ à valeurs complexes,
nous rappelons la définition de la covariance de $f$ et $g$~:
$$\mbox{Cov}(f,g)=
       {\bf E}\left[fg\right]-{\bf E}\left[f\right]
       {\bf E}\left[g\right].$$
Nous considérons la suite de variables aléatoires $(X_k)_{k\ge 0}$
définie par $X_k=f\circ T^k$ où $f:\Omega\rightarrow\bf R$
est une fonction mesurable bornée et $\nu$-centrée.

Pour tout entier naturel $n$, nous notons $S_n:=\sum_{k=1}^nX_k$, 
avec la convention $S_0=0$.
Rappelons tout d'abord le résultat suivant (prouvé dans \cite{cjan1}) 
qui résulte de la preuve de Rio \cite{Rio}.
\begin{thm1}[\cite{Rio}]\label{sec:thm0}
Soit $(X_k)_{k\ge 0}$ une suite stationnaire 
de variables aléatoires réelles bornées
centrées 
définies sur un même espace probabilisé telle que, 
pour tous entiers naturels $a,b,c$ vérifiant $a+b+c\le 3$
on a~:
\begin{equation}\label{sec:a1}
 \sum_{j\ge 1}j.\sup_{i\le i+j\le i+k\le i+l}
  \sup_{F\in\cal L}\left\Vert {\bf E}\left[{X_{i+j}}^{a}{X_{i+k}}^{b}
  {X_{i+l}}^{c}\vert X_0,...,X_i  \right]-{\bf E}
   \left[{X_{i+j}}^{a}{X_{i+k}}^{b}
  {X_{i+l}}^{c}\right]\right\Vert_\infty<+\infty.
\end{equation}
Alors, la limite suivante existe~:
$$\sigma:=\lim_{n\rightarrow +\infty} {1\over \sqrt{n}}
       \left({\bf E}\left[{S_n}^2\right]\right)^{1\over 2}.$$
Si $\sigma=0$, alors la suite $(S_n)_n$ est bornée dans $L^2$.

Si $\sigma>0$, alors la suite de variables aléatoires
$\left({S_n\over\sqrt{n}}\right)_{n\ge 1}$ converge en loi vers 
une variable aléatoire $N$ de loi
normale centrée de variance $\sigma^2$ et il existe un nombre
réel $B>0$ tel que, pour tout entier $n\ge 1$, on ait~:
$$\sup_{x\in\bf R}\left\vert \nu\left({S_n\over\sqrt n}\le x\right)
   -{\bf P}(N\le x)\right\vert \le {B\over \sqrt{n}}.$$
\end{thm1}
En étudiant en détail la preuve de Rio, on s'aperçoit qu'il utilise en fait 
la propriété suivante (qui découle de la propriété (\ref{sec:a1})
et peut la remplacer dans l'énoncé de son théorème)~:
pour toute fonction continue $F:{\bf R}^3\rightarrow\bf R$
et tous entiers $0\le j\le k\le l\le l+p\le l+q\le l+s$, on a~:
$$\left\vert Cov\left(F\left(S_{j-1},X_j,X_k,X_l\right),
    {X_{l+p}}^{a}{X_{l+q}}^{b}{X_{l+s}}^{c}\right)
    \right\vert\le
     C\left\Vert F\left(S_{j-1},X_j,X_k,X_l\right)
 \right\Vert_{L^1}\varphi_p\ \mbox{avec}\ \sum_{p\ge 1}p\varphi_p<+\infty.$$
En fait, on peut encore affaiblir les hypothèses.
En utilisant la méthode de Rio, nous allons montrer
le résultat suivant.
\begin{thm1}\label{sec:thmZ}
Soit $(X_k)_{k\ge 0}$ une suite stationnaire de variables 
aléatoires réelles bornées
centrées 
définies sur un même espace probabilisé telle qu'il
existe trois nombres réels $C\ge 1$, $M
\ge\max\left(1,\Vert X_0\Vert_\infty\right)$ 
et $r\ge 0$ 
et une suite de nombres réels 
$(\varphi_{p,l})_p$
majorées par 1 vérifiant 
$\sum_{p\ge 1}p\max_{l=0,...,{p\over r+1}}\varphi_{p,l}<+\infty$
tels que pour tous entiers naturels $a,b,c$ vérifiant $a+b+c\le 3$, 
pour tous entiers
$j,k,l,p,s$ vérifiant~: $ 1\le j\le k\le l\le l+p\le l+q\le l+s$,
pour toute fonction différentiable $F:{\bf R}^3\rightarrow\bf R$
telle que $F\left(S_{j-1},X_j,X_k,X_l\right)$ soit intégrable,
nous avons~:

\noindent$\displaystyle \left\vert Cov\left(F\left(S_{j-1},X_j,X_k,X_l\right),
    {X_{l+p}}^{a}{X_{l+q}}^{b}{X_{l+s}}^{c}\right)
    \right\vert\le$
\begin{equation}\label{sec:a2}  
\le C\left(\left\Vert F\left(S_{j-1},X_j,X_k,X_l\right)
       \right\Vert_{L^1}+
       \left\Vert\sup_{\vert u\vert,\vert v\vert,\vert w\vert,\vert x\vert\le M}
        \vert DF\left(S_{j-1}+u,X_j+v,X_k+w,X_l+x\right)\vert_\infty
        \right\Vert_{L^1}\right)
        \varphi_{p,s-p} .
\end{equation}
Alors, on a les mêmes conclusions que dans le théorème \ref{sec:thm0}.
\end{thm1}

%
Remarquons tout d'abord que, la suite $(X_k)_{k\ge 0}$ étant stationnaire,
pour tout $n\ge 1$, on a~:
\begin{equation}\label{sec:covariance}
{\bf E}\left[\left({S_n\over\sqrt{n}}\right)^2\right]=
   {\bf E}[{X_0}^2]+2\sum_{k=1}^{n}\left(1-{ k\over n}\right)
    {\bf E}[X_0X_k]. 
\end{equation}
La condition de décorrélation (\ref{sec:a2}) assure donc l'existence de $\sigma$
et le fait que l'on a~:
$$\sigma^2={\bf E}\left[{X_0}^2\right]+2\sum_{k\ge 1}{\bf E}\left[
      X_0X_k\right]. $$
De plus, d'après la formule (\ref{sec:covariance}), on a~:
$$ \left\vert{\bf E}\left[{S_n}^2\right]
  -n\sigma^2\right\vert\le 2\sum_{k\ge 1}k\left\vert {\bf E}[X_0X_k]
    \right\vert\le 4CM\sum_{k\ge 1}k\varphi_{k,0}.$$
Ainsi, si $\sigma=0$, alors la suite de variables aléatoires
$(S_n)_{n\ge 1}$ est bornée dans $L^2(\Omega,{\bf R})$.

Supposons à présent $\sigma>0$. Dans la suite, nous supposons
que $\sigma=1$ (on peut toujours s'y ramener, quitte à remplacer $X_k$
par $X_k\over\sigma$).

\subsection{Les grandes lignes de la preuve}

L'idée de la preuve est d'utiliser un raisonnement par récurrence.
Avant de nous lancer dans la preuve technique du théorème,
nous en expliquons la démarche.

Nous voulons montrer l'existence d'un nombre réel $A\ge 1$ tel que la propriété
$\left({\cal P}_n(A)\right)$ 
suivante soit vérifiée pour tout entier $n\ge 2$~:
$$\left({\cal P}_n(A)\right)\ : \ \forall k=1,...n-1,\ \sup_{x\in{\bf R}}
\left\vert \nu\left({S_k\over\sqrt k}\le x\right)
   -{\bf P}(N\le x)\right\vert \le {A\over\sqrt{k}}.$$

Posons $V_n:={\bf E}\left[{S_n}^2\right]$ et
$v_n:={\bf E}\left[{S_n}^2\right]-{\bf E}\left[{S_{n-1}}^2\right]$.
Nous avons alors~:
\begin{eqnarray*}
v_n&=&{\bf E}\left[(S_{n-1}+X_n)^2-{S_{n-1}}^2\right]\\
&=&{\bf E}\left[X_n\left(2S_{n-1}+X_n\right)\right]\\
&=&{\bf E}[{X_0}^2]+2\sum_{k=1,...,n-1}{\bf E}\left[X_0X_k\right].
\end{eqnarray*}
Ainsi la suite $(v_n)_{n\ge 1}$ converge vers $\sigma^2=1$.
Il existe donc un entier $n_0\ge 1$ tel que, pour tout $n\ge n_0$, on a~:
$v_n\ge {1\over 2}$.
Dans la suite, nous supposons l'existence d'une suite $(N_i)_{i\ge 0}$ 
de variables aléatoires indépendantes identiquement distribuées 
de loi normale centrée réduite et indépendante de la suite de 
variables aléatoires $(X_k)_{k\ge 0}$.
L'essentiel de la preuve est l'établissement du résultat suivant~:
\begin{prop}\label{sec:clef}
Sous les hypothèses du théorème \ref{sec:thmZ},
 il existe un nombre réel $K\ge M$ 
et une fonction $\psi:[1,+\infty]\rightarrow ]0;+\infty[$ 
continue décroissante
vérifiant $\lim_{\varepsilon\rightarrow +\infty}\psi(\varepsilon)=0$ tels
que pour tout entier $n\ge 6n_0$ et tout nombre réel $A\ge M$, si on a
$\left({\cal P}_n(A)\right)$, alors on a, pour tout $\varepsilon\ge 1$,~:
$$\sup_{x\in\bf R}\left\vert{\nu}\left(\{S_n+\varepsilon Y\le x\}\right)
     -{\nu}\left(\{S_{n_0-1}+T_{n_0-1,n}+\varepsilon Y\le x\}\right)
     \right\vert \le {K\over\sqrt{n}}\left(1+A\psi(\varepsilon)\right),$$
où $Y$ et $T_{n_0-1,n}$ sont deux variables aléatoires indépendantes
et indépendantes de $(X_k)_{k\ge 0}$ de lois normales centrées de variances
respectives $1$ et $V_n-V_{n_0}$.
\end{prop}
Montrons à présent comment ce résultat nous permet de conclure.
\begin{enumerate}
\item 
Nous remarquons tout d'abord que, pour tout $A\ge \sqrt{6n_0}$,
la propriété $\left({\cal P}_{6n_0}(A)\right)$ est vérifiée.
\item 
Montrons qu'il existe un nombre
réel $A_0\ge M$ tel que, pour tout entier $n\ge 6n_0$ et tout 
nombre réel $A\ge A_0$, on a~:
$\left({\cal P}_{n}(A)\right)\Rightarrow \left({\cal P}_{n+1}(A)\right)$.

Soit un entier $n\ge 6n_0$ et un nombre réel $A\ge M$ tels qu'on ait
$\left({\cal P}_n(A)\right)$.
Alors, d'après la proposition \ref{sec:clef}, nous avons~:
$$\sup_{x\in\bf R}\left\vert{\nu}\left(\{S_n+\varepsilon Y\le x\}\right)
     -{\nu}\left(\{S_{n_0-1}+T_{n_0-1,n}+\varepsilon Y\le x\}\right)
     \right\vert \le {K\over\sqrt{n}}\left(1+A\psi(\varepsilon)\right).$$
D'autre part, si $T_n$ est une variable aléatoire de loi normale 
${\cal N}(0,V_n)$ indépendante de $Y$, alors nous avons~:
$$ \left\vert {\nu}\left(\{S_{n_0-1}+T_{n_0-1,n}+
  \varepsilon Y\le x\}\right)-{\nu}\left(\{T_n+
  \varepsilon Y\le x\}\right)\right\vert\le C_0{n_0 M}^2{1\over n}.$$
Ainsi, il existe une constante $K'\ge1$ (indépendante de $n\ge 6n_0$ et 
de $A\ge M$) telle que~:
$$\sup_{x\in\bf R}\left\vert{\nu}\left(\{S_n+\varepsilon Y\le x\}\right)
     -{\nu}\left(\{T_n+\varepsilon Y\le x\}\right)
     \right\vert \le {K'\over\sqrt{n}}\left(1+A\psi(\varepsilon)\right).$$
Notons $\varepsilon_A$ l'unique réel $\varepsilon_A\in\left[{
   K'\over 3},+\infty\right[$ tel que ${7\over 2}
   \left(K'(1+A\psi(\varepsilon_A))\right)
   =3\varepsilon_A$.
L'existence et l'unicité de $\varepsilon_A$ résulte
de la décroissance de la fonction $\varepsilon\mapsto
{7\over 2}K'(1+A\psi(\varepsilon))$ et de la stricte croissance de 
$\varepsilon\mapsto 3\varepsilon$ sur $\left[{
   K'\over 3},+\infty\right[$.
D'après le lemme 2  de Rio, nous en déduisons que nous avons~:
$$\sup_{x\in\bf R}\left\vert{\nu}\left(\{S_n\le x\}\right)
     -{\nu}\left(\{T_n\le x\}\right)
     \right\vert \le {3\varepsilon_A\over\sqrt{n}}=
       {7\over 2} {K'\over\sqrt{n}}\left(1+A\psi(\varepsilon_A)\right.$$
Nous avons ainsi montré l'existence d'un nombre réel $K"\ge 1$ tel que,
pour tout entier $n\ge 6n_0$ et tout nombre réel $A\ge M$, nous avons~:
$$\left({\cal P}_n(A)\right)\ \ \Rightarrow\ \ 
\left({\cal P}_{n+1}\left(K"(1+A\psi(\varepsilon_A))\right)\right),$$
où $\varepsilon_A$ est l'unique réel $\varepsilon_A\in\left[{
   K\over 3},+\infty\right[$ tel que $
   \left(K"(1+A\psi(\varepsilon_A))\right)
   =3\varepsilon_A$.

Puis nous montrons qu'il existe un $A_0$ tel que, pour tout $A\ge A_0$,
on a $3\varepsilon_A\le A$.
En effet, on a $\lim_{A\rightarrow +\infty}\varepsilon_A=+\infty$
car $A\mapsto\varepsilon_A$ est croissante et
ne peut-être majorée sinon $A\mapsto\psi(\varepsilon_A)$
serait minorée par un nombre $m>0$, et on aurait~:
$3\varepsilon_A\ge K"A\psi(\varepsilon_A)\ge K"Am$.
Il existe donc un $A_0$ tel que, pour tout $A\ge A_0$, on a~:
$K"\le{A\over 2}$ et $K"A\psi(\varepsilon_A)\le{A\over 2}$
et donc $3\varepsilon_A=K"(1+A\psi(\varepsilon_A))\le A$.
\end{enumerate}
Ainsi, pour tout nombre réel $A\ge\max(A_0,\sqrt{6n_0})$, on a~:
$\left({\cal P}_{6n_0}(A)\right)$ et, pour tout $n\ge 6n_0$,
$\left({\cal P}_{n}(A)\right)\Rightarrow \left({\cal P}_{n+1}(A)\right)$.

\subsection{La démonstration de la proposition \ref{sec:clef}}

Soit un entier $n\ge n_0$.
Nous nous donnons une suite de variables aléatoires indépendantes 
$(Y_k)_{k\ge n_0}$ définies sur $(\Omega,{\cal F},\nu)$
indépendante de $(X_k)_{k\ge 0}$ telle que 
$Y_k$ est de loi normale centrée de variance $v_k$.
Considérons également une variable
aléatoire $Y$ de loi normale centrée réduite
définie sur $(\Omega,{\cal F},\nu)$ indépendante
de $\left((Y_k)_{k\ge n_0},(X_k)_{k\ge 0}\right)$
(on peut toujours le faire quitte à se placer dans un espace plus gros).
\begin{nota}
On pose~: $\Delta_k(f)={\bf E}\left[f(S_{k-1}+X_k)\right]
-{\bf E}\left[f(S_{k-1}+Y_k)\right]$.

Pour tous nombres réels $\varepsilon>0$ et $y$,
on pose~: $f_{k,\varepsilon,y}(x)={\nu}\left(\varepsilon Y+\sum_{i=k+1}^nY_i+
x\le y\right),$
avec la convention usuelle $\sum_{i=k+1}^nY_i=0$ si $k=n$.
\end{nota}
Remarquons que l'on a~:
$$ {\nu}\left(S_n+\varepsilon Y\le y\right)-{\nu}
\left(S_{n_0-1}+\varepsilon Y+\sum_{i=n_0}^{n}Y_i\le y\right)
   =\sum_{k=n_0}^n\Delta_{k}\left(f_{k,\varepsilon,y}\right).$$

Nous écrivons~: $\Delta_k(f)=\Delta_{1,k}(f)-\Delta_{2,k}(f)$,
avec~:
$$\Delta_{1,k}(f)={\bf E}\left[f(S_{k-1}+X_k)\right]
-{\bf E}\left[f(S_{k-1})\right]-{v_k\over 2}{\bf E}[f"(S_{k-1})] $$
et
$$\Delta_{2,k}(f)={\bf E}\left[f(S_{k-1}+Y_k)\right]
-{\bf E}\left[f(S_{k-1})\right]-{v_k\over 2}{\bf E}[f"(S_{k-1})].$$
Nous présentons à présent les grandes lignes de la preuve de la proposition
\ref{sec:clef}. Puis nous démontrerons en détail chacun des
résultats intermédiaires énoncés ci-dessous
(à l'exception du lemme \ref{sec:lemme7}).
\subsubsection{Quand $k$ n'est pas grand}

\begin{lem}[adaptation du lemme 6 de \cite{Rio}]\label{sec:lemme6}
Soit une fonction $f:{\bf R}\rightarrow\bf R$ de classe $C^3$. Nous
avons~:
$$\left\vert\Delta_k(f)\right\vert\le 
\Vert f'''\Vert_\infty \left(M^3+20
  C^2M^{2}(r+1)\sum_{p=0}^{k-1}(1+p)\varphi_{p,0}
  +3CM\sum_{p=r+1}^{k-2}
   \sum_{l=1,...,k-1 : (r+1)l\le p\ }\varphi_{p,l}\right).$$
\end{lem}
\begin{lem}[lemme 5 de \cite{Rio}]\label{sec:lemme5}
Pour tout entier $k\ge n_0$,
la fonction $f_{k,\varepsilon,y}$ est infiniment dérivable et, pour tout
entier $i\ge 1$, nous avons~:
$$\Vert {f_{k,\varepsilon,y}}
  ^{(i)}\Vert_\infty\le {C_i\over (n-k+\varepsilon^2)^
   {i\over 2}} ,$$
en notant $C_i:=2^{i\over 2}\Vert \phi^{(i-1)}\Vert_\infty$, où 
$\phi$ désigne la fonction densité de la loi normale centrée réduite.
\end{lem}

\begin{prop}\label{sec:prop0}
Pour tous nombres réels $y$ et $\varepsilon\ge 1$, on a~:
$$\sum_{k=n_0}^{n-\left\lfloor{n\over 3}\right\rfloor-1}
\left\vert\Delta_{k}\left(f_{k,\varepsilon,y}\right)\right\vert\le
   {1\over\sqrt{n}}{8\sqrt{3}\over\sqrt{\pi e}}
 \left(M^3+20
  C^2M^{2}\sum_{p=0}^{k-1}(1+p)\varphi_{p,0}
+3CM\sum_{p=r+1}^{k-2}
   \sum_{l=1}^{\left\lfloor{p\over r+1}\right\rfloor}\varphi_{p,l}\right).$$
\end{prop}

{\bf Preuve de la proposition \ref{sec:prop0}.\/}
Remarquons tout d'abord que, d'après le lemme \ref{sec:lemme5}, nous avons~:
$$\Vert {f_{k,\varepsilon,y}}'''\Vert_\infty\le {4\over\sqrt{\pi e}}
{1\over (n-k+\varepsilon^2)^{3\over 2}}. $$
Ainsi, d'après le lemme \ref{sec:lemme6}, on a~:
\begin{eqnarray*}
\sum_{k=n_0}^{n-\left\lfloor{n\over 3}\right\rfloor-1}
\left\vert\Delta_k(f_{k,\varepsilon,y})\right\vert&\le& 
{4\over\sqrt{\pi e}}\sum_{k=\left\lfloor {n\over 3}\right\rfloor+1}^{n-1}
{1\over (k+1)^{3\over 2}} \left(M^3+20
  C^2M^{2}\sum_{p=0}^{k-1}(1+p)\varphi_{p,0}\right.\\
&\ &\ \ \ \ \ \ \ \ \ \ \ \ \ \ \ \ \ 
  \left.+3CM\sum_{p=r+1}^{k-2}
   \sum_{l=1,...,k-1 : (r+1)l\le p\ }\varphi_{p,l}\right)\\
&\le&{1\over\sqrt{n}}{8\sqrt{3}\over\sqrt{\pi e}}
 \left(M^3+20
  C^2M^{2}\sum_{p=0}^{k-1}(1+p)\varphi_{p,0}\right.\\
&\ &\ \ \ \ \ \ \ \ \ \ \ \ \ \ \ \ \ 
  \left.+3CM\sum_{p=r+1}^{k-2}
   \sum_{l=1,...,k-1 : (r+1)l\le p\ }\varphi_{p,l}\right),
\end{eqnarray*}
{\bf cqfd.}
\subsubsection{Quand les $k$ sont grands}

\noindent Nous allons à présent nous intéresser aux quantités
$\Delta_k\left(f_{k,\varepsilon,y}\right)$ lorsque $ k\ge n-
  \left\lfloor{n\over 3}\right\rfloor$.
Nous allons utiliser les mêmes idées que celles de la preuve du lemme
\ref{sec:lemme6} (formules de Taylor et sommes glissantes).
Nous aurons besoin de majorations plus fines (faisant intervenir
des contrôles en norme 1 au lieu de contrôles en norme infinie).

La preuve du lemme suivant est identique à celle du lemme 7 de \cite{Rio}.
Elle utilise le lemme \ref{sec:lemme5} et les majorations 
établies dans les lemmes 3 et 4 de \cite{Rio} sous des hypothèses générales.
\begin{lem}[adaptation du lemme 7 de \cite{Rio}]\label{sec:lemme7}
Pour tout entier $i\ge 1$,
il existe une constante $K_i$ telle que, pour tout entier $n\ge 6n_0$,
tout entier  
$k\ge n-\left\lfloor{n\over 3}\right\rfloor$, tout entier
$l\in\left[{n\over 3};k\right]$, tout nombre réel $A\ge M$, si la propriété
$\left({\cal P}_n(A)\right)$ est vérifiée, alors on a~: 
\begin{equation}\label{sec:majo1}
\sup_{a\in\bf R}{\bf E}\left[\sup_{z\in[-3M;3M]}
\left\vert{f_{k,\varepsilon,y}}^{(i)}
    (S_l+a-z)\right\vert\right]\le {K_i\over\sqrt{n}}
    \left({A\over (n-k+\varepsilon^2)^{i\over 2}}+
    {1\over (n-k+\varepsilon^2)^{i-1\over 2}}\right)
\end{equation}
et
\begin{equation}\label{sec:majo2}
\sup_{a\in\bf R}\left\vert{\bf E}\left[
{f_{k,\varepsilon,y}}^{(i)}
    (S_l+a)\right]\right\vert\le {K_i\over\sqrt{n}}
    \left({A\over (n-k+\varepsilon^2)^{i\over 2}}+
    {1\over n^{i-1\over 2}}\right).
\end{equation}
\end{lem}

Remarquons tout d'abord que, d'après une estimation
de $\left\vert\Delta_{2,k}(f)\right\vert$ établie dans la
preuve du lemme \ref{sec:lemme6} (cf. formule (\ref{sec:d2k}) et le lemme
\ref{sec:lemme7}, nous avons le résultat suivant.
\begin{lem}\label{sec:lemme00}
Soit un entier $n\ge 6n_0$ et un nombre réel $A\ge M$.
Si la propriété $\left({\cal P}_n(A)\right)$ est vérifiée, alors,
pour tout entier $k=n-\left\lfloor {n\over 3}\right\rfloor,...,n-1$, on a~:
\begin{equation}
\left\vert\Delta_{2,k}\left(f_{k,\varepsilon,y}\right)\right\vert\le
{4\sqrt{\pi}\over 3\sqrt{6}}
  {K_3\over\sqrt{n}}
    \left({A\over (n-k+\varepsilon^2)^{3\over 2}}+
    {1\over n}\right)(CM)^{3\over 2}
    \sum_{p=0}^{k-1}(1+p)\varphi_{p,0}.
\end{equation}
\end{lem}

Intéressons nous à présent à $\Delta_{1,k}\left(f_{
k,\varepsilon,y}\right)$. Le contrôle de cette quantité est plus difficile.
Nous établissons le résultat suivant.
\begin{lem}\label{sec:lemme000}
Il existe un nombre réel $\tilde K$ (ne dépendant que de 
$K_1$, $K_2$, $K_3$, $K_4$, $C$, $M$ et $r$) tel que,
pour tout entier $n\ge 6n_0$ et tout nombre réel $A\ge M$,
si la propriété $\left({\cal P}_n(A)\right)$ est vérifiée, alors,
pour tout entier $k=n-\left\lfloor {n\over 3}\right\rfloor,...,n-1$, on a~:
\begin{equation}
\left\vert\Delta_{1,k}\left(f_{k,\varepsilon,y}\right)\right\vert\le
{\tilde K\over\sqrt{n}}\left({A\alpha_{ \sqrt{n-k}}
  \over (n-k+\varepsilon^2)^{3\over 2}}
 +{\beta_{ \sqrt{n-k}}\over n}
 +{A\gamma_{\sqrt{n-k}}\over n-k+\varepsilon^2}
 +{A\delta_{\sqrt{n-k}}\over\sqrt{n-k+\varepsilon^2}}
    +\varphi_{\lfloor \sqrt{n-k}\rfloor+1,0}\right),
\end{equation}
avec $\alpha_m:=1+\sum_{p= 1}^{\lfloor m\rfloor}p\zeta_p$,
$\beta_m:=1+\sum_{p= 1}^{\lfloor m\rfloor}\zeta_p$ et
$\gamma_m:=\sum_{p=\left\lfloor{m\over (r+2)^2}\right\rfloor}
   ^{\lfloor m\rfloor}\zeta_p$ et 
   $\delta_{m}:=\sum_{p=\lfloor m\rfloor+1}^{+\infty}
      {\zeta_p\over p}$, en notant $\zeta_p:=p\max_
      {j=0,...,\left\lfloor{p\over r+1}\right\rfloor}\varphi_{p,j}$.
\end{lem}

\subsubsection{Conclusion}
Achevons la preuve de la proposition \ref{sec:clef}.
Soit un entier $n\ge 6n_0$ et un nombre réel $A\ge M$
tels que $\left({\cal P}_n(A)\right)$ soit vérifiée.
Soit un nombre réel $\varepsilon\ge 1$ et un nombre réel $y$
quelconque.

Rappelons que nous avons~:
$${\nu}\left(S_n+\varepsilon Y\le y\right)-
  {\nu}\left(S_{n_0-1}+T_{n_0-1,n}+\varepsilon Y\le y\right)=
  \sum_{k=n_0}^{n}
 \Delta_{k}\left(f_{k,\varepsilon,y}\right),$$
en notant $T_{n_0-1,n}:=\sum_{i=n_0}^nY_i$.
Alors, d'après la proposition \ref{sec:prop0}, nous avons~:
$$\sum_{k=n_0}^{n-\left\lfloor {n\over 3}\right\rfloor-1}
\left\vert \Delta_{k}\left(f_{k,\varepsilon,y}\right)\right\vert
\le {K"_1\over\sqrt{n}},$$
où la constante $K"_1$ est indépendante de $(n,A,\varepsilon,y)$.

D'autre part, d'après le lemme \ref{sec:lemme00}, on a~:
\begin{equation}
\sum_{k=n-\left\lfloor{n\over 3}\right\rfloor}
  ^{n}\left\vert\Delta_{2,k}\left(f_{k,\varepsilon,y}\right)\right\vert
 \le 
  {K"_2\over\sqrt{n}}\left(1+A\sum_{l\ge 0}
      {1\over (l+\varepsilon^2)^{3\over 2}}\right).
\end{equation}
D'autre part, d'après le lemme \ref{sec:lemme000}, nous avons~:
$$
\left\vert\Delta_{1,k}\left(f_{k,\varepsilon,y}\right)\right\vert\le
{\tilde K\over\sqrt{n}}\left({A\alpha_{ \sqrt{n-k}}
  \over (n-k+\varepsilon^2)^{3\over 2}}
 +{\beta_{ \sqrt{n-k}}\over n}
 +{A\gamma_{\sqrt{n-k}}\over n-k+\varepsilon^2}
 +{A\delta_{\sqrt{n-k}}\over\sqrt{n-k+\varepsilon^2}}
    +\varphi_{\lfloor \sqrt{n-k}\rfloor+1,0}\right),
$$
avec $\alpha_m=1+\sum_{p= 1}^{\lfloor m\rfloor}p\zeta_p$,
$\beta_m=1+\sum_{p= 1}^{\lfloor m\rfloor}\zeta_p$ et
$\gamma_m=\sum_{p=\left\lceil{m\over (r+1)^2}\right\rceil}
   ^{\lfloor m\rfloor}\zeta_p$ et 
   $\delta_{m}=\sum_{j=\lfloor m\rfloor+1}^{+\infty}
      {\zeta_{j}\over j}$, en notant $\zeta_p:=p\max_
      {j=0,...,\left\lfloor{p\over r+1}\right\rfloor}\varphi_{p,j}$.
Nous faisons la somme sur $k\in\left\{n-\left\lfloor {n\over 3}\right\rfloor
,...,n\right\}$ et contrôlons chacun des termes.
\begin{enumerate}
\item Contrôle du premier terme~:
\begin{eqnarray*}
\sum_{k=n-\left\lfloor{n\over 3}\right\rfloor}^n 
{A\alpha_{ \sqrt{n-k}}
  \over (n-k+\varepsilon^2)^{3\over 2}}&=&A\sum_{l=0}^{\left\lfloor{n\over 3}\right\rfloor} 
{1+\sum_{p= 1}^{\lfloor \sqrt{l}\rfloor}p\zeta_p
  \over (l+\varepsilon^2)^{3\over 2}}\\
&\le& A\left(\left(\sum_{l\ge 0}
{1 \over (l+\varepsilon^2)^{3\over 2}}\right)+
  \sum_{p=1}^{\sqrt{n}}\left(
   \sum_{l\ge p^2}
    {1\over (l+\varepsilon^2)^{3\over 2}}\right)p\zeta_p\right).
\end{eqnarray*}
D'où~:
\begin{equation}
\sum_{k=n-\left\lfloor{n\over 3}\right\rfloor}^n 
{A\alpha_{ \sqrt{n-k}}
  \over (n-k+\varepsilon^2)^{3\over 2}}\le A\left(\left(\sum_{l\ge 0}
{1 \over (l+\varepsilon^2)^{3\over 2}}\right)+
  \sum_{p=1}^{\sqrt{n}}
    {2\over \sqrt{p^2+\varepsilon^2-1}}p\zeta_p\right).
\end{equation}
\item Contrôle du second terme~:
\begin{equation}
\sum_{k=n-\left\lfloor{n\over 3}\right\rfloor}^n 
{\beta_{\sqrt{n-k}}\over n}
\le 1+\sum_{p\ge 1}\zeta_p.
\end{equation}
\item Contrôle du troisième terme~:
\begin{eqnarray*}
\sum_{k=n-\left\lfloor{n\over 3}\right\rfloor}^n{A\gamma_{\sqrt{n-k}}
  \over n-k+\varepsilon^2}&\le&
  \sum_{l=0}^{\left\lfloor{n\over 3}\right\rfloor}{A\gamma_{\sqrt{l}}
  \over l+\varepsilon^2}\\
  &\le&
 A\sum_{p=0}^{\left\lfloor\sqrt{n}\right\rfloor}
    \sum_{l=p^2}^{p^2(r+2)^4}{1\over l+\varepsilon^2}\zeta_p
\end{eqnarray*}
D'où~:
\begin{equation}
\sum_{k=n-\left\lfloor{n\over 3}\right\rfloor}^n{A\gamma_{\sqrt{n-k}}
  \over n-k+\varepsilon^2}
  \le A\sum_{p\ge 0}\ln\left({p^2(r+2)^4+\varepsilon^2
     \over p^2-1+\varepsilon^2}\right)\zeta_p.
\end{equation}
\item Contrôle du quatrième terme~:
\begin{eqnarray*}
\sum_{k=n-\left\lfloor{n\over 3}\right\rfloor}^n{A\delta_{\sqrt{n-k}}
  \over \sqrt{n-k+\varepsilon^2}}&\le&A
  \sum_{l=0}^{\left\lfloor{n\over 3}\right\rfloor}{1
  \over \sqrt{l+\varepsilon^2}}
  \sum_{p=\lfloor\sqrt{l}\rfloor+1}^{+\infty}
      {\zeta_p\over p}\\
      &\le& A\sum_{p\ge 1}\sum_{l=0}^{p^2}{1
  \over \sqrt{l+\varepsilon^2}}{\zeta_p\over p}\\
  &\le& A\sum_{p\ge 1}2\left(\sqrt{p^2+\varepsilon^2}
    -\sqrt{\varepsilon^2-1}\right){\zeta_p\over p}.
\end{eqnarray*}
Ainsi, on a~:
\begin{equation}
\sum_{k=n-\left\lfloor{n\over 3}\right\rfloor}^n{A\delta_{\sqrt{n-k}}
  \over\sqrt{ n-k+\varepsilon^2}}
 \le A\sum_{p\ge 1}{2(1+p^2)\over\sqrt{p^2+\varepsilon^2}}
   {\zeta_p\over p}.
\end{equation}
\item Contrôle du dernier terme~:
\begin{eqnarray*}
\sum_{k=n-\left\lfloor{n\over 3}\right\rfloor}^n
    \varphi_{\lfloor \sqrt{n-k}\rfloor+1,0}&=&\sum_{l=0}
    ^{\left\lfloor{n\over 3}\right\rfloor}
    \varphi_{\lfloor \sqrt{l}\rfloor+1,0}\\
&\le&\sum_{p\ge 1}\#\left\{l : \lfloor\sqrt{l}\rfloor=p\right\}
     \varphi_{p,0}.
\end{eqnarray*}
D'où~:
\begin{equation}\sum_{k=n-\left\lfloor{n\over 3}\right\rfloor}^n
    \varphi_{\lfloor \sqrt{n-k}\rfloor+1,0}\le\sum_{p\ge 1}(2p+1)
     \varphi_{p,0}.
\end{equation}
\end{enumerate}
{\it cqfd.}
\section{Les calculs pour les petits $k$}
Rappelons l'énoncé du lemme \ref{sec:lemme6}~:

\noindent{\bf Lemme \ref{sec:lemme6} } {\it 
Soit une fonction $f:{\bf R}\rightarrow\bf R$ de classe $C^3$. Nous
avons~:
$$\left\vert\Delta_k(f)\right\vert\le 
\Vert f'''\Vert_\infty \left(M^3+20
  C^2M^{2}(r+1)\sum_{p=0}^{k-1}(1+p)\varphi_{p,0}
  +3CM\sum_{p=r+1}^{k-2}
   \sum_{l=1,...,k-1 : (r+1)l\le p\ }\varphi_{p,l}\right).$$
}\medskip

\noindent{\it Preuve du lemme \ref{sec:lemme6}.\/}
Nous écrivons à nouveau~: $\Delta_k(f)=\Delta_{1,k}(f)-\Delta_{2,k}(f)$,
avec~:
$$\Delta_{1,k}(f)={\bf E}\left[f(S_{k-1}+X_k)\right)
-{\bf E}\left[f(S_{k-1})\right]-{v_k\over 2}{\bf E}[f"(S_{k-1})] $$
et
$$\Delta_{2,k}(f)={\bf E}\left[f(S_{k-1}+Y_k)\right)
-{\bf E}\left[f(S_{k-1})\right]-{v_k\over 2}{\bf E}[f"(S_{k-1})].$$
Comme $Y_k$ est indépendante de $S_{k-1}$ et est de loi normale
centrée de variance $v_k$, on a~:
\begin{eqnarray*}
\left\vert\Delta_{2,k}(f)\right\vert&=&\left\vert
    {\bf E}\left[f\left(S_{k-1}+Y_k\right)
   -f\left(S_{k-1}\right)-{1\over 2}f"\left(S_{k-1}\right){Y_k}^2\right]
     \right\vert\\
 &\le&\left\vert{\bf E}\left[{1\over 2}\int_0^1(1-t)^2
    f'''(S_{k-1}+tY_k)\, dt. {Y_k}^3
           \right]\right\vert\\
 &\le&{1\over 2}\int_0^1(1-t)^2
   \sup_{a\in\bf R}\left\vert{\bf E}\left[ 
       f'''(S_{k-1}+a)\right]\right\vert{\bf E}\left[\vert{Y_k}\vert
       ^3\right] \, dt \\
 &\le&{1\over 6}\sup_{a\in\bf R}\left\vert{\bf E}\left[ 
       f'''(S_{k-1}+a)\right]\right\vert
   {\bf E}\left[\vert{Y_k}\vert
       ^3\right]  \\
 &\le& {2\over 3\sqrt{2\pi}}\sup_{a\in\bf R}\left\vert{\bf E}\left[ 
       f'''(S_{k-1}+a)\right]\right\vert{v_k}^{3\over 2},
\end{eqnarray*}
car ${\bf E}\left[\left\vert Y_k\right\vert^3\right]\le{4{v_k}^{3\over 2}
     \over\sqrt{2\pi}}$.
Or, d'après l'expression de $v_k$ et l'hypothèse (\ref{sec:a2}), on a~:
$v_k\le 2CM\sum_{p=0}^{k-1}\varphi_{p,0}$.
En utilisant l'inégalité de Hölder et le fait que $\varphi_{p,0}\le 1$,
on peut montrer (cf. \cite{Rio} page 264) que l'on a~:
$${v_k}^{3\over 2}\le (4CM)^{3\over 2}
    {\pi\over \sqrt{6}}\sum_{p=0}^{k-1}(1+p)\varphi_{p,0}.$$
Ainsi, nous avons~:
\begin{equation}\label{sec:d2k}
\left\vert\Delta_{2,k}(f)\right\vert\le
  {8\sqrt{\pi}\over 3\sqrt{3}}\sup_{a\in\bf R}\left\vert{\bf E}\left[ 
       f'''(S_{k-1}+a)\right]\right\vert
  (CM)^{3\over 2}\sum_{p=0}^{k-1}(1+p)\varphi_{p,0}.
\end{equation}
Nous allons à présent contrôler $\Delta_{1,k}(f)$.
Nous avons~:
\begin{eqnarray*}
\Delta_{1,k}(f)&=&{\bf E}\left[f\left(S_{k-1}+X_k\right)
   -f\left(S_{k-1}\right)-{1\over 2}f"\left(S_{k-1}\right)v_k\right]\\
&=&{\bf E}\left[f'\left(S_{k-1}\right)X_{k}
  +{1\over 2}f"\left(S_{k-1}\right)\left({X_k}^2-v_k\right)+{1\over 6}f'''\left(S_{k-1}+\theta_kX_k\right){X_k}^3\right],
\end{eqnarray*}
où 
$\theta_k$ est une variable aléatoire à valeurs dans $[0;1]$.
Rappelons d'autre part que~:
$$v_k-{\bf E}\left[{X_k}^2\right]=
   2\sum_{i=1}^{k-1}{\bf E}\left[X_iX_k\right].$$
D'où~:
\begin{eqnarray*}
\Delta_{1,k}(f)&=&
{\bf E}\left[f'(S_{k-1})X_k\right]+
{1\over 2} Cov\left(f"(S_{k-1}),{X_k}^2\right) \\
&-&{\bf E}[f"(S_{k-1})]
       \sum_{i=1}^{k-1}{\bf E}[X_iX_k]\\
&+& {\bf E}\left[{1\over 6}f'''\left(S_{k-1}+\theta_kX_k\right)
     {X_k}^3\right].\
\end{eqnarray*}
Nous allons contrôler un à un chacun des termes du membre de 
droite de cette égalité.
On a~:
\begin{equation}
\left\Vert{1\over 6}f'''\left(S_{k-1}+\theta_kX_k\right){X_k}^3 
\right\Vert_\infty\le{1\over 6}\Vert f'''\Vert_\infty M^3.
\end{equation}
D'autre part on a~:
\begin{eqnarray*}
\left\vert Cov\left(f"(S_{k-1}),{X_k}^2\right)\right\vert&\le&\sum_{i=1}^{k-1}\left\vert
    Cov\left(f"(S_i)-f"(S_{i-1}),{X_k}^2\right)\right\vert\\
&\le& C\left(\Vert f'''\Vert_\infty M+2\Vert f'''\Vert_\infty\right)
     \sum_{p=1}^{k-1}\varphi_{p,0},
\end{eqnarray*}
en utilisant (\ref{sec:a2}).
D'où~:
\begin{equation} \left\vert Cov\left(f"(S_{k-1}),{X_k}^2\right)    \right\vert\le 3CM\Vert f'''\Vert_\infty
     \sum_{p=1}^{k-1}\varphi_{p,0}.
\end{equation}
Il nous reste à contrôler la quantité ${\bf E}\left[f'(S_{k-1})X_k\right]-{\bf E}[f"(S_{k-1})]
       \sum_{i=1}^{k-1}{\bf E}[X_iX_k] $.
On a~:
$$f'(S_{k-1})=\sum_{i=1}^{k-1}\left(f'(S_i)-f'(S_{i-1})\right)=
    \sum_{i=1}^{k-1}\left(f"(S_{i-1})X_i+\int_{S_{i-1}}^{S_i}(S_i-t)f'''(t)\, dt\right). $$
Nous avons donc~:
\begin{equation}\label{sec:A1}
{\bf E}\left[f'(S_{k-1})X_k\right]-{\bf E}[f"(S_{k-1})]
       \sum_{i=1}^{k-1}{\bf E}[X_iX_k] =
       \sum_{i=1}^{k-1}\left(Cov\left(f"(S_{i-1}),X_iX_k\right)\right.
\end{equation}
       $$+\sum_{i=1}^{k-1}{\bf E}\left[f"(S_{i-1})-f"(S_{k-1})\right]{\bf E}[X_iX_k]$$
   $$ \left.+\sum_{i=1}^{k-1}Cov\left(\int_{S_{i-1}}^{S_i}(S_i-t)
   f'''(t)\, dt,X_k\right) \right).$$
\`A l'aide de (\ref{sec:a2}), nous obtenons~:
$$\left\vert{\bf E}\left[f"(S_{i-1})-f"(S_{k-1})\right]
    \right\vert\left\vert{\bf E}[X_iX_k]\right\vert\le 
    \Vert f'''\Vert_\infty C M^2(k-i)
     \varphi_{k-i,0}$$
et donc~:
$$
\left\vert\sum_{i=1}^{k-1}{\bf E}[X_iX_k] {\bf E}
[f"(S_{i-1})-f"(S_{k-1})]\right\vert\leq
\Vert f'''\Vert_\infty C M^2\sum_{p=1}^{k-1} p\varphi_{p,0}.
$$
De plus, à l'aide de (\ref{sec:a2}), nous avons~:
$$\left\vert Cov\left(\int_{S_{i-1}}^{S_i}(S_i-t)f'''(t)\, dt,X_k
  \right)\right\vert\le
   C\left(\Vert f'''\Vert_\infty M^2+4\Vert f'''\Vert_\infty M\right)
     \varphi_{k-i,0}$$
et donc~:
\begin{equation}
\sum_{i=1}^{k-1}\left\vert Cov\left(\int_{S_{i-1}}^{S_i}(S_i-t)
     f'''(t)\, dt,X_k\right)\right\vert\le
   5C\Vert f'''\Vert_\infty M^2\sum_{p=1}^{k-1}\varphi_{p,0}.
\end{equation}
Il nous reste à contrôler le premier morceau du membre de droite
de (\ref{sec:A1}).
Pour tout entier $i=1,...,k-1$, en notant 
$j=j_i:=\max(0,(r+2)i-(r+1)k)$, à l'aide de \ref{sec:a2}, nous
avons~:
\begin{eqnarray*}
\left\vert Cov\left(\left(f"(S_{i-1})-f"(S_j))X_i,X_k\right)\right)\right\vert 
&\le& \sum_{m=j+1}^{i-1}\left\vert Cov\left((f"(S_{m})-f"(S_{m-1})
      )X_i,X_k\right)\right\vert\\
&\le&5C\Vert f'''\Vert_\infty M^2(i-j-1)\varphi_{k-i,0}
\end{eqnarray*}
et
$$\left\vert{\bf E}\left[f"(S_{i-1})-f"(S_j)\right]
    \right\vert{\bf E}[X_iX_k]\le \Vert f'''\Vert_\infty C M^2(i-j-1)
     \varphi_{k-i,0} .$$
Ainsi, on a~:
\begin{equation}
\sum_{i=1}^{k-1}\left\vert Cov\left((f"(S_{i-1})-f"(S_j)),
     X_iX_k\right)\right\vert\le 6CM^2\Vert f'''\Vert_\infty
       (r+1)\sum_{p=1}^{k-1}
                p\varphi_{p,0}.
\end{equation}
Si $(r+2)i-(r+1)k\le 0$, alors $j=0$ et donc $Cov\left(f"(S_j),X_iX_k\right)=0$.
Pour tout $i=1,...,k$ tel que $(r+2)i-(r+1)k>0$, on a~:
\begin{eqnarray*}
\left\vert Cov\left(f"(S_j),X_iX_k\right)\right\vert&\le&
 \sum_{l=1}^j\left\vert Cov\left(f"(S_l)-f"(S_{l-1}),X_iX_k\right)
 \right\vert\\
 &\le& 3CM\Vert f'''\Vert_\infty\sum_{l=1}^{(r+2)i-(r+1)k}
 \varphi_{i-l,k-i}\\
 &\le&3CM\Vert f'''\Vert_\infty\sum_{p=(r+1)(k-i)}^{i-1}\varphi_{p,k-i}.
\end{eqnarray*}
Ainsi, on a~:
\begin{eqnarray}
\sum_{i=1}^{k-1}\left\vert Cov\left(f"(S_j),X_iX_k\right)\right\vert&\le&
3CM\Vert f'''\Vert_\infty\sum_{i=1,...,k-1:(r+1)k<(r+2)i\ }
   \sum_{p=(r+1)(k-i)}^{i-1}\varphi_{p,k-i}\nonumber\\
&\le& 3CM\Vert f'''\Vert_\infty\sum_{j=1,...,k-1:(r+1)j<k\ }
   \sum_{p=(r+1)j}^{k-l-1}\varphi_{p,j}\nonumber\\
&\le& 3CM\Vert f'''\Vert_\infty\sum_{p=r+1}^{k-2}
   \sum_{j=1,...,k-1 : (r+1)j\le p\ }\varphi_{p,j}.
\end{eqnarray}
{\it cqfd.}\medskip

\section{Les calculs pour les grands $k$}

Nous voulons prouver le lemme \ref{sec:lemme000} dont nous rappelons l'énoncé~:
\medskip

\noindent{\bf Lemme \ref{sec:lemme000} } {\it
Il existe un nombre réel $\tilde K$ (ne dépendant que de 
$K_1$, $K_2$, $K_3$, $K_4$, $C$, $M$ et $r$) tel que,
pour tout entier $n\ge 6n_0$ et tout nombre réel $A\ge M$,
si la propriété $\left({\cal P}_n(A)\right)$ est vérifiée, alors,
pour tout entier $k=n-\left\lfloor {n\over 3}\right\rfloor,...,n-1$, on a~:
$$
\left\vert\Delta_{1,k}\left(f_{k,\varepsilon,y}\right)\right\vert\le
{\tilde K\over\sqrt{n}}\left({A\alpha_{ \sqrt{n-k}}
  \over (n-k+\varepsilon^2)^{3\over 2}}
 +{\beta_{ \sqrt{n-k}}\over n}
 +{A\gamma_{\sqrt{n-k}}\over n-k+\varepsilon^2}
 +{A\delta_{n,k}\over\sqrt{n-k+\varepsilon^2}}
    +\varphi_{\lfloor \sqrt{n-k}\rfloor+1,0}\right),
$$
avec $\alpha_m:=1+\sum_{p= 1}^{\lfloor m\rfloor}p\zeta_p$,
$\beta_m:=1+\sum_{p= 1}^{\lfloor m\rfloor}\zeta_p$ et
$\gamma_m:=\sum_{p=\left\lfloor{m\over (r+2)^2}\right\rfloor}
   ^{\lfloor m\rfloor}\zeta_p$ et 
   $\delta_{n,k}:=\sum_{p=\lfloor\sqrt{n-k}\rfloor+1}^k
      {\zeta_p\over p}$, en notant $\zeta_p:=p\max_
      {j=0,...,\left\lfloor{p\over r+1}\right\rfloor}\varphi_{p,j}$.
}\medskip

Nous ferons un usage répété du résultat suivant dont la preuve est identique 
à celle du lemme 7 de \cite{Rio}.
\medskip

\noindent{\bf Lemme \ref{sec:lemme7} } {\it 
Pour tout entier $i\ge 1$,
il existe une constante $K_i$ telle que, pour tout entier $n\ge 6n_0$,
tout entier  
$k\ge n-\left\lfloor{n\over 3}\right\rfloor$, tout entier
$l\in\left[{n\over 3};k\right]$, tout nombre réel $A\ge M$, si la propriété
$\left({\cal P}_n(A)\right)$ est vérifiée, alors on a~: 
$$
\sup_{a\in\bf R}{\bf E}\left[\sup_{z\in[-3M;3M]}
\left\vert{f_{k,\varepsilon,y}}^{(i)}
    (S_l+a-z)\right\vert\right]\le {K_i\over\sqrt{n}}
    \left({A\over (n-k+\varepsilon^2)^{i\over 2}}+
    {1\over (n-k+\varepsilon^2)^{i-1\over 2}}\right)
$$
et
$$
\sup_{a\in\bf R}\left\vert{\bf E}\left[
{f_{k,\varepsilon,y}}^{(i)}
    (S_l+a)\right]\right\vert\le {K_i\over\sqrt{n}}
    \left({A\over (n-k+\varepsilon^2)^{i\over 2}}+
    {1\over n^{i-1\over 2}}\right).
$$
}\medskip

\noindent{\it Preuve du lemme \ref{sec:lemme000}.\/}
Par souci de lisibilité, nous
noterons dans la suite $f_k$ pour $f_{k,\varepsilon,y}$.
Commençons par écrire $\Delta_{1,k}(f_k)$ comme une somme de termes 
que nous majorerons ensuite séparément.
Un développement limité avec reste intégral donne :
\begin{eqnarray*}
\Delta_{1,k}(f_k)&=&{\bf E}\left[f_k(S_k)-f_k(S_{k-1})-{v_k\over 2}
     f"(S_{k-1})\right]\\
  &=&{\bf E}\left[{f_k}'(S_{k-1})X_k\right]+{1\over 2}{\bf E}
     \left[{f_k}"(S_{k-1})({X_k}^2-v_k)\right]+{1\over 6}
     {\bf E}
     \left[{f_k}'''(S_{k-1}){X_k}^3\right]\\
&{ }& \ \ \ \ +{\bf E}\left[{{X_k}^4\over 3!}\int_0^1(1-t)^3{f_k}^{(4)}
\left(S_{k-1}
+tX_k\right)\, dt\right].\\
\end{eqnarray*}
Comme le nombre $v_k$ vaut ${\bf E}[{X_0}^2]+2\sum_{i=1,...,k-1}
{\bf E}\left[X_0X_i\right]$, nous avons 
$$
{1\over 2}{\bf E}
     \left[{f_k}"(S_{k-1})({X_k}^2-v_k)\right]={1\over 2} 
     Cov({f_k}"(S_{k-1}),X_k^2)-\sum_{i=1,...,k-1}
     {\bf E}\left[{f_k}"(S_{k-1})\right]{\bf E}
     \left[X_kX_i\right].
$$
Nous avons donc :
\begin{eqnarray*}
\Delta_{1,k}(f_k)&=&{\bf E}\left[{f_k}'(S_{k-1})X_k\right]-
\sum_{i=1,...,k-1}{\bf E}\left[X_kX_i\right]{\bf E}\left[{f_k}"(S_{k-1})
    \right]\\
&+&{1\over 2} Cov({f_k}"(S_{k-1}),X_k^2)\\
&+&{1\over 6}
     {\bf E}
     \left[{f_k}'''(S_{k-1}){X_k}^3\right]\\
&+& {\bf E}\left[{{X_k}^4\over 3!}\int_0^1(1-t)^3f^{(4)}\left(S_{k-1}
+tX_k\right)\, dt\right].
\end{eqnarray*}
Pour tout entier naturel $l$ vérifiant $l\leq k-1$, on peut écrire :
\begin{eqnarray*}
{f_k}'(S_{k-1})-{f_k}'(S_{k-l-1})
&=&\sum_{j=k-l}^{k-1}{f_k}'(S_{j})-{f_k}'(S_{j-1})\\
&=& \sum_{j=k-l}^{k-1}{f_k}"(S_{j-1})X_j\\
&+& {1\over 2}\sum_{j=k-l}^{k-1}{f_k}'''(S_{j-1})X_j^2\\
&+& {1\over 2}\sum_{j=k-l}^{k-1}{X_{j}}^3\int_0^1(1-t)^2{f_k}^{(4)}
        \left(S_{j-1}+tX_{j}\right)\, dt.
\end{eqnarray*}
En remplaçant ${f_k}'(S_{k-1})$ par l'expression donnée par 
cette égalité on obtient :
\begin{eqnarray*}
{\bf E}\left[{f_k}'(S_{k-1})X_k\right]-\sum_{i=1,...,k-1}
{\bf E}\left[X_kX_i\right]{\bf E}\left[{f_k}"(S_{k-1})\right]
&=& {\bf E}\left[{f_k}'(S_{k-l-1})X_k\right]\\
&+& \sum_{j=k-l}^{k-1}{\bf E}\left[{f_k}"(S_{j-1})
\left(X_jX_k-{\bf E}\left[X_kX_j\right]\right)\right]\\
&+& \sum_{j=k-l}^{k-1}
  {\bf E}\left[{f_k}"(S_{k-1})-{f_k}"(S_{j-1})\right]
      {\bf E}\left[X_kX_j\right]\\
&+&  \sum_{i=1}^{k-l-1}{\bf E}[{f_k}"(S_{k-1})]{\bf E}\left[X_kX_i\right]\\
&+& {1\over 2}\sum_{j=k-l}^{k-1}{\bf E}\left[{f_k}'''(S_{j-1})X_j^2X_k\right]\\
&+& {1\over 2}\sum_{j=k-l}^{k-1}{\bf E}\left[X_k{X_{j}}^3
    \int_0^1(1-t)^2{f_k}^{(4)}
        \left(S_{j-1}+tX_{j}\right)\, dt\right].\\
\end{eqnarray*}
Nous avons ainsi exprimé $\Delta_{1,k}(f_k)$ en une somme 
\begin{eqnarray*}
\Delta_{1,k}(f_k)
&=& {\bf E}\left[{{X_k}^4\over 3!}\int_0^1(1-t)^3f^{(4)}\left(S_{k-1}
+tX_k\right)\, dt\right]\\
&+&{1\over 6}
     {\bf E}
     \left[{f_k}'''(S_{k-1}){X_k}^3\right]\\
&+& 
{\bf E}\left[{f_k}'(S_{k-l-1})X_k\right]\\
&+&  \sum_{i=1}^{k-l-1}{\bf E}\left[{f_k}"(S_{k-1})\right]
    {\bf E}\left[X_kX_i\right]\\
&+&{1\over 2} Cov({f_k}"(S_{k-1}),X_k^2)\\
&+& {1\over 2}\sum_{j=k-l}^{k-1}{\bf E}\left[{f_k}'''(S_{j-1})X_j^2X_k\right]\\
&+& {1\over 2}\sum_{j=k-l}^{k-1}{\bf E}\left[X_k{X_{j}}^3
     \int_0^1(1-t)^2{f_k}^{(4)}
        \left(S_{j-1}+tX_{j}\right)\, dt\right]\\
&+& \sum_{j=k-l}^{k-1}{\bf E}\left[{f_k}"(S_{k-1})-{f_k}"(S_{j-1})\right]
     {\bf E}\left[X_kX_j\right]\\
&+& \sum_{j=k-l}^{k-1}{\bf E}\left[{f_k}"(S_{j-1})\left(X_jX_k-
    {\bf E}\left[X_kX_j\right]\right)\right].
\end{eqnarray*}
Nous allons maintenant majorer, dans l'ordre, chacun des termes de cette somme.
{\bf Dorénavant l'entier $l$ est pris égal à $\lfloor\sqrt{n-k}\rfloor$ la partie entière de $\sqrt{n-k}$.}

\subsection{Contrôle de ${\bf E}[{{X_k}^4\over 3!}\int_0^1(1-t)^3f^{(4)}\left(S_{k-1}
+tX_k\right)\, dt]$}

D'après la formule (\ref{sec:majo1}) du lemme \ref{sec:lemme7}, nous
avons~:
\begin{eqnarray*}
\left\vert {\bf E}\left[{{X_k}^4\over 3!}\int_0^1(1-t)^3f^{(4)}\left(S_{k-1}
+tX_k\right)\, dt\right]\right\vert&\le&
  {1\over 24}{\bf E}\left[\sup_{u\in[-M;M]}\left\vert
       {f_k}^{(4)}(S_{k-1}+u) \right\vert\right]M^4\\
     &\le& {K_4 M^4\over 24\sqrt{n}}
    {2A\over (n-k+\varepsilon^2)^{3\over 2}}.
\end{eqnarray*}
Nous avons donc~:
\begin{equation}
\left\vert {\bf E}\left[{{X_k}^4\over 3!}\int_0^1(1-t)^3f^{(4)}\left(S_{k-1}
+tX_k\right)\, dt\right]\right\vert\le {K_4 M^4\over 12\sqrt{n}}
    {A\alpha_{\sqrt{n-k}}\over (n-k+\varepsilon^2)^{3\over 2}}.
\end{equation}

\subsection{Contrôle de ${1\over 6}
     {\bf E}\left[{f_k}'''(S_{k-1}){X_k}^3\right]$} 
Nous avons~:
$${f_k}'''(S_{k-1})={f_k}'''(S_{k-l-1})+\sum_{j=1}^l\left({f_k}'''(S_{k-j})
     -{f_k}'''(S_{k-j-1})\right) .$$
D'après (\ref{sec:a2}) et (\ref{sec:majo1}), nous avons~:
\begin{eqnarray}
\left\vert Cov\left({f_k}'''(S_{k-l-1}),{X_k}^3\right)\right\vert
&\le&C\left({2K_3A\over\sqrt{n}(n-k+\varepsilon^2)}
    +{2K_4A\over\sqrt{n}(n-k+\varepsilon^2)^{3\over 2}}\right)\varphi_{l+1,0}
       \nonumber\\
&\le& {2C(K_3+K_4)\over\sqrt{n}}{A\delta_{{n,k}}\over \sqrt
  {n-k+\varepsilon^2}}\label{sec:B1}
\end{eqnarray}
et
\begin{eqnarray}
\sum_{j=1}^l\left\vert Cov\left(\left({f_k}'''(S_{k-j})
     -{f_k}'''(S_{k-j-1})\right),{X_k}^3)\right]\right\vert&\le&
\sum_{j=1}^l{6K_4M\over\sqrt n}{A\over (n-k+\varepsilon^2)^{3\over2}}
\varphi_{j,0}\nonumber\\
&\le&{6K_4M\over\sqrt n}
{A\alpha_{\sqrt{n-k}}\over (n-k+\varepsilon^2)^{3\over2}}.\label{sec:B2}
\end{eqnarray}
De plus, d'après (\ref{sec:majo2}), nous avons~:
\begin{eqnarray}
\left\vert {\bf E}\left[{f_k}'''(S_{k-1})
     \right]{\bf E}\left[{X_k}^3\right]\right\vert
&\le&{K_3\over\sqrt n}\left({A\over (n-k+\varepsilon^2)^{3\over 2}}
  +{1\over n}\right)M^3\nonumber\\
&\le&{K_3M^3\over\sqrt n}\left({A\alpha_{\sqrt{n-k}}
  \over (n-k+\varepsilon^2)^{3\over 2}}
  +{\beta_{\sqrt{n-k}}\over n}\right)\label{sec:B3}.
\end{eqnarray}

\subsection{Contrôle de ${\bf E}\left[{f_k}'(S_{k-l-1})X_k\right]$}
\`A l'aide de (\ref{sec:a2}) et (\ref{sec:majo1}), nous avons~:
\begin{eqnarray}
\left\vert {\bf E}\left[{f_k}'(S_{k-l-1})X_k\right]\right\vert
&=&\left\vert Cov\left({f_k}'(S_{k-l-1}),X_k\right)\right\vert\nonumber\\
&\le& {C(K_1+2K_2)\over\sqrt{n}}
\left({A\over \sqrt{n-k+\varepsilon^2}}+1\right)
\varphi_{l+1,0}\nonumber\\
&\le&{C(K_1+2K_2)\over\sqrt{n}}\left({A\delta_{n,k}\over \sqrt{n-k+\varepsilon^2}}
   +\varphi_{\lfloor\sqrt{n-k}\rfloor+1,0}\right)\label{sec:B4}.
\end{eqnarray}

\subsection{Contrôle de $
  \sum_{j=l+1}^{k-1}\left\vert 
    {\bf E}\left[{f_k}''(S_{k-1})\right]{\bf E}[X_{k-j}X_k]\right\vert$}
  
Nous avons~:
\begin{eqnarray}
 \sum_{j=l+1}^{k-1}
  \left\vert{\bf E}\left[{f_k}''\left(S_{k-1}\right)
    \right]{\bf E}\left[X_{k-j}X_k\right]\right\vert
    &=& \sum_{j=l+1}^{k-1}\left\vert {\bf E}\left[{f_k}''\left(S_{k-1}\right)
    \right]Cov\left(X_{k-j},X_k\right)\right\vert\nonumber\\
    &\le& \sum_{j=l+1}^{k-1}{2CK_2\over\sqrt{n}}{A\over\sqrt{n-k+\varepsilon^2}
      }2CM\varphi_{j,0}\nonumber\\
    &\le& {4CMK_2\over\sqrt{n}}{A\delta_{n,k}\over\sqrt{n-k+\varepsilon^2}}
       .\label{sec:B5}
\end{eqnarray}

\subsection{Contrôle de $\left\vert Cov\left(
{f_k}''(S_{k-1}),{X_k}^2\right)\right\vert$}
Nous avons :
\begin{eqnarray*}
{f_k}''(S_{k-1}) &=&{f_k}''(S_{k-l-1})+\sum_{j=1}^l
\left({f_k}''(S_{k-j})-{f_k}''(S_{k-j-1})\right)\\
&=&{f_k}''(S_{k-l-1})+\sum_{j=1}^l
{f_k}'''(S_{k-j-1})X_{k-j}\\
&{ }& \ \ \ +\sum_{j=1}^l{X_{k-j}}^2\int_0^1(1-t){f_k}^{(4)}
        \left(S_{k-j-1}+tX_{k-j}\right)\, dt.\\
\end{eqnarray*}
Ainsi, on a~:
\begin{eqnarray*}
Cov\left({f_k}''(S_{k-1}),{X_k}^2\right)&=&Cov\left(
{f_k}''(S_{k-l-1}),{X_k}^2\right)+\sum_{j=1}^lCov\left(
{f_k}'''(S_{k-j-1})X_{k-j},{X_k}^2\right)\\
&+&\sum_{j=1}^l\int_0^1(1-t)Cov\left(
{f_k}^{(4)}(S_{k-j-1}+tX_{k-j}){X_{k-j}}^2,{X_k}^2\right)\, dt
\end{eqnarray*}
Nous allons à présent contrôler chaque terme du membre
de droite de cette égalité.
\begin{enumerate}
\item Tout d'abord, d'après (\ref{sec:a2}) et (\ref{sec:majo1}),
nous avons~:
\begin{eqnarray*}
\left\vert Cov\left(
{f_k}''(S_{k-l-1}),{X_k}^2\right)\right\vert&\le&
C{2(K_2+K_3)\over\sqrt{n}}{A\over\sqrt{n-k+\varepsilon^2}}\varphi_{l+1,0}\\
&\le& C{2(K_2+K_3)\over\sqrt{n}}
{A\delta_{n,k}\over\sqrt{n-k+\varepsilon^2}}.
\end{eqnarray*}
\item Nous contrôlons à présent la quantité $\sum_{j=1}^l 
Cov\left(
{f_k}'''(S_{k-j-1})X_{k-j},{X_k}^2\right)$.

Pour tout $j=1,...,l$, on a~:
$${f_k}'''(S_{k-j-1})={f_k}'''(S_{k-l-1})+
   \sum_{m=j+1}^l\left({f_k}'''(S_{k-m})-
      {f_k}'''(S_{k-m-1})\right).$$
\begin{enumerate}
\item Pour tout couple d'entiers $(j,m)$ vérifiant $1\le j<m\le l$
et $m\le (r+2)j$, d'après (\ref{sec:a2}) et (\ref{sec:majo1}),
nous avons~:
$$\left\vert Cov\left(\left({f_k}'''(S_{k-m})-
      {f_k}'''(S_{k-m-1})\right)X_{k-j},{X_k}^2\right)\right\vert
      \le {18 CK_4 M^2\over{\sqrt{n}}}{A\over (n-k+\varepsilon^2)^{3\over 2}}
       \varphi_{j,0}$$
Ainsi, nous avons~:
\begin{equation}
\sum_{j=1}^l\sum_{m=j+1}^{\min\left(l,(r+2)j\right)}
 \left\vert Cov\left(\left({f_k}'''(S_{k-m})-
      {f_k}'''(S_{k-m-1})\right)X_{k-j},{X_k}^2\right)\right\vert
      \le {18 CK_4 M^2(r+1)\over{\sqrt{n}}}{A\alpha_{\sqrt{n-k}}
       \over (n-k+\varepsilon^2)^{3\over 2}}.
\end{equation}
\item Nous souhaitons à présent majorer la quantité suivante~:
$$\left\vert\sum_{j=1}^l\sum_{m=\min\left(l,(r+2)j\right)+1}^l
  Cov\left(\left({f_k}'''(S_{k-m})-
      {f_k}'''(S_{k-m-1})\right)X_{k-j},{X_k}^2\right)\right\vert.$$
Remarquons que cette quantité est en fait~:
$$\left\vert\sum_{j=1}^l
  Cov\left(\left({f_k}'''(S_{k-\min\left(l,(r+2)j\right)-1})-
      {f_k}'''(S_{k-l-1})\right)X_{k-j},{X_k}^2\right)\right\vert.$$
Nous allons utiliser la formule~: 
$Cov(AB,C)=Cov(A,BC)-Cov(A,B){\bf E}[C]+{\bf E}[A]Cov(B,C)
   $. 
Pour tout couple d'entiers $(j,m)$ vérifiant $1\le j\le(r+2)j+1\le 
m\le l$, d'après (\ref{sec:a2}) et (\ref{sec:majo1}),
nous avons~:
$$\left\vert Cov\left({f_k}'''(S_{k-m})-
      {f_k}'''(S_{k-m-1}),X_{k-j}{X_k}^2\right)\right\vert\le
C{6K_4M\over\sqrt n}{A\over(n-k+\varepsilon^2)^{3\over2}}\varphi_{m-j,j}
$$
et
$$\left\vert Cov\left({f_k}'''(S_{k-m})-
      {f_k}'''(S_{k-m-1}),X_{k-j}\right)
      {\bf E}\left[{X_k}^2\right]\right\vert
       \le C{6K_4M\over\sqrt{n}}{A\over(n-k+\varepsilon^2)^{3\over2}}
         \varphi_{m-j,0}.$$
D'où nous obtenons~:
\begin{equation}
\sum_{j=1}^l\sum_{(r+2)j+1\le m\le l}\left\vert Cov\left({f_k}'''(S_{k-m})-
      {f_k}'''(S_{k-m-1}),X_{k-j}{X_k}^2\right)\right\vert\le
{6CK_4M\over (r+1)\sqrt n}{A\alpha_{\sqrt{n-k}}
    \over(n-k+\varepsilon^2)^{3\over2}} 
\end{equation}
et
\begin{equation}
\sum_{j=1}^l\sum_{(r+2)j+1\le m\le l}
\left\vert Cov\left({f_k}'''(S_{k-m})-
      {f_k}'''(S_{k-m-1}),X_{k-j}\right)
      {\bf E}\left[{X_k}^2\right]\right\vert
       \le {6CK_4M\over(r+1)\sqrt{n}}{A\alpha_{\sqrt{n-k}}
       \over(n-k+\varepsilon^2)^{3\over2}}.
\end{equation}
En effet, nous avons~:
\begin{eqnarray*}
\sum_{j=1}^l\sum_{(r+1)j+1\le m\le l}\varphi_{m-j,j}&=&
\sum_{j=1}^{\left\lfloor{l\over r+2}\right\rfloor}
\sum_{p=(r+1)j+1}^{l-j}
   \varphi_{p,j}\\
&=&\sum_{p=r+2}^l\sum_{j=1}^{\left\lfloor {p\over r+1}\right\rfloor}
  \varphi_{p,j}\\
&\le& {\alpha_{\sqrt{n-k}}\over r+1}
\end{eqnarray*}
et, de même~:
$$ \sum_{j=1}^l\sum_{(r+2)j+1\le m\le l}\varphi_{m-j,0}\le 
{\alpha_{\sqrt{n-k}}\over r+1}.$$
Enfin, en utilisant (\ref{sec:a2}) et (\ref{sec:majo2}), nous obtenons~:

$\displaystyle
\sum_{j=1}^{l}\left\vert{\bf E}\left[{f_k}'''(S_{k-\min\left(l,(r+2)j
    \right)-1})-
      {f_k}'''(S_{k-l-1})\right] Cov(X_{k-j},{X_k}^2)\right\vert$
\begin{eqnarray}
&=&\sum_{j=1}^{\left\lfloor{l\over r+2}\right\rfloor}
\left\vert{\bf E}\left[{f_k}'''(S_{k-(r+2)j
    -1})-
      {f_k}'''(S_{k-l-1})\right] Cov(X_{k-j},{X_k}^2)\right\vert
      \nonumber \\
&\le& \sum_{j=1}^{\left\lfloor{l\over r+2}\right\rfloor}
 {2K_3\over\sqrt{n}}\left({A\over(n-k+\varepsilon^2)^{3\over 2}}
    +{1\over n}\right)2CM\varphi_{j,0}\nonumber\\
&\le& {4CK_3M\over\sqrt{n}}\left({A\alpha_{\sqrt{n-k}}
   \over(n-k+\varepsilon^2)^{3\over 2}}
    +{\beta_{\sqrt{n-k}}\over n}\right)
\end{eqnarray}
\item Pour tout entier $j=1,...,l$ vérifiant ${l\over{r+2}}< j\le l$,
on a~: $\sqrt{n-k}<(r+2)j$ et,
d'après (\ref{sec:a2}) et (\ref{sec:majo1}), nous avons~:
$$\left\vert Cov\left(  {f_k}'''(S_{k-l-1})X_{k-j},{X_k}^2
  \right)\right\vert\le {4C(K_3+K_4)M\over\sqrt{n}}
   {A\over{n-k+\varepsilon^2}}\varphi_{j,0}.$$
D'où
\begin{equation}
\sum_{j=\left\lfloor{\sqrt{n-k}\over r+2}\right\rfloor+1}^l
\left\vert Cov\left(  {f_k}'''(S_{k-l-1})X_{k-j},{X_k}^2
  \right)\right\vert\le  {4C(K_3+K_4)M\over\sqrt{n}}
     {A\gamma_{\sqrt{n-k}}\over{n-k+\varepsilon^2}}.
\end{equation}
\item Pour tout entier $j=1,...,l$ vérifiant $j\le {l\over{r+2}}$,
on a~: $j\le{\sqrt{n-k}\over r+2}$ et, 
d'après (\ref{sec:a2}) et (\ref{sec:majo1}), 
nous avons~:
$$\left\vert Cov\left(  {f_k}'''(S_{k-l-1}),X_{k-j}{X_k}^2
  \right)\right\vert\le {2C(K_3+K_4)\over\sqrt{n}}
   {A\over{n-k+\varepsilon^2}}\varphi_{l+1-j,0} $$
et
$$\left\vert Cov\left(  {f_k}'''(S_{k-l-1}),X_{k-j})\right)
     {\bf E}\left[{X_k}^2
  \right]\right\vert\le {2CM^2(K_3+K_4)\over\sqrt{n}}
   {A\over{n-k+\varepsilon^2}}\varphi_{l+1-j,0}. $$
De plus, d'après (\ref{sec:a2}) et (\ref{sec:majo2}), nous avons~:
$$\left\vert {\bf E}\left[  {f_k}'''(S_{k-l-1})\right]
{\bf E}\left[X_{k-j}{X_k}^2
  \right]\right\vert\le {K_3\over\sqrt{n}}
   \left({A\over{(n-k+\varepsilon^2)^{3\over 2}}}+{1\over n}\right)
   C2M\varphi_{j,0}. $$
D'où
\begin{equation}
\sum_{j=1}^{\left\lfloor {\sqrt{n-k}\over{r+2}}\right\rfloor}
\left\vert Cov\left(  {f_k}'''(S_{k-l-1})X_{k-j},{X_k}^2
  \right)\right\vert\le {4CM^2(K_3+K_4)\over\sqrt{n}}
   {A\gamma_{\sqrt{n-k}}\over{n-k+\varepsilon^2}}+
\end{equation}
$$+ {2CMK_3\over\sqrt{n}}
   {A\alpha_{\sqrt{n-k}}\over{(n-k+\varepsilon^2)^{3\over 2}}}
   +{K_3\over\sqrt{n}}{\beta_{\sqrt{n-k}}\over n}$$
\end{enumerate}
\item Nous avons~:

$\displaystyle\sum_{j=1}^l\int_0^1(1-t)\left\vert Cov\left(
{f_k}^{(4)}(S_{k-j-1}+tX_{k-j})X_{k-j}^2,{X_k}^2\right)\right\vert\, dt\le$
\begin{eqnarray}
&\le& \sum_{j=1}^l\int_0^1(1-t) {10(K_4+K_5)M^2\over\sqrt{n}}
{A\over{(n-k+\varepsilon^2)}^{3\over 2}}\varphi_{j,0}\, dt \nonumber\\
&\le& {5(K_4+K_5)M^2\over\sqrt{n}}{A\alpha_{\sqrt{n-k}}\over
{(n-k+\varepsilon^2)}^{3\over 2}}.\label{sec:B6}
\end{eqnarray}
\end{enumerate}

\subsection{Contrôle de $\sum_{j=1}^l{\bf E}\left[{f_k}'''(S_{k-j-1}){X_{k-j}}^2X_k\right]$}
Nous contrôlons cette quantité de 
la même manière que nous avons contrôlé 
$$
\sum_{j=1}^lCov\left({f_k}'''(S_{k-j-1})X_{k-j},{X_k}^2\right)
$$
dans la section précédente.
Nous obtenons ainsi des majorations analogues (avec des 
constantes sensiblement différentes).

\subsection{Contrôle de $\sum_{j=1}^l\int_0^1(1-t)
{\bf E}\left[{f_k}^{(4)}(S_{k-j-1}+tX_{k-j}){X_{k-j}}^3X_k\right]\, dt$}

Nous contrôlons la quantité 
$\sum_{j=1}^l\int_0^1(1-t)
Cov\left({f_k}^{(4)}
(S_{k-j-1}+tX_{k-j}){X_{k-j}}^3,X_k\right)\, dt$
comme nous avons contrôlé
$\sum_{j=1}^l\int_0^1(1-t)
Cov\left({f_k}^{(4)}(S_{k-j-1}){X_{k-j}}^2,{X_k}^2\right)$.
Nous obtenons ainsi une majoration analogue à (\ref{sec:B6}).
\subsection{Contrôle de $\sum_{j=1}^l{\bf E}\left[{f_k}''(S_{k-j-1})
   -{f_k}''(S_{k-1})\right]{\bf E}\left[X_{k-j}X_k\right]$}

Pour tout entier $j=1,...,l$, on a~:
\begin{eqnarray*}
{f_k}''(S_{k-1})
   -{f_k}''(S_{k-j-1})&=&\sum_{m=1}^j
  \left({f_k}''(S_{k-m})
   -{f_k}''(S_{k-m-1})\right)\\
&=&\sum_{m=1}^j
  {f_k}'''(S_{k-m-1})X_{k-m}+X_{k-m}^2\int_0^1(1-t)f^{(4)}\left(S_{k-m-1}
+tX_{k-m}\right)\, dt.
\end{eqnarray*}
Nous utilisons l'égalité 
$${f_k}'''(S_{k-m-1})={f_k}'''(S_{k-l-1})+\sum_{p=m+1}^l
     \left({f_k}'''(S_{k-p})-{f_k}'''(S_{k-p-1})\right) $$
pour aboutir à l'expression :
\begin{eqnarray*}
&{ }&\sum_{j=1}^l{\bf E}\left[{f_k}''(S_{k-j-1})
   -{f_k}''(S_{k-1})\right]{\bf E}\left[X_{k-j}X_k\right]\\
&{ }&\ \ \ \ \ \ \ \ =\sum_{j=1}^l\sum_{m=1}^j
   {\bf E}\left[X_{k-m}^2\int_0^1(1-t)f^{(4)}\left(S_{k-m-1}
+tX_{k-m}\right)\, dt\right]{\bf E}\left[X_{k-j}X_k\right]\\
&{ }&\ \ \ \ \ \ \ \ +\sum_{j=1}^l\sum_{m=1}^j{\bf E}\left[{f_k}'''(S_{k-l-1})X_{k-m}\right]
     {\bf E}\left[X_{k-j}X_k\right]\\
&{ }&\ \ \ \ \ \ \ \ +\sum_{j=1}^l\sum_{m=1}^j\sum_{p=m+1}^{l}
   {\bf E}\left[\left({f_k}'''(S_{k-p})-{f_k}'''(S_{k-p-1})\right)
   X_{k-m}\right]{\bf E}\left[X_{k-j}X_k\right].\\
\end{eqnarray*}
\begin{enumerate}
\item En utilisant (\ref{sec:a2}) et (\ref{sec:majo1}), nous obtenons~:

$\displaystyle\sum_{j=1}^l\sum_{m=1}^j\left\vert
   {\bf E}\left[X_{k-m}^2\int_0^1(1-t)f^{(4)}\left(S_{k-m-1}
+tX_{k-m}\right)\, dt\right]{\bf E}\left[X_{k-j}X_k\right]\right\vert\le$
\begin{eqnarray}
&\le&\sum_{j=1}^l\sum_{m=1}^j{K_4\over\sqrt{n}}{A\over(n-k+\varepsilon^2)
   ^{3\over 2}}M^2CM\varphi_{j,0}  \nonumber \\
&\le&
   {CM^3K_4\over\sqrt{n}}{A\over (n-k+\varepsilon^2)^{3\over 2}}
     \alpha_{\sqrt{n-k}}.
\end{eqnarray}
\item
Contrôlons la somme :

$\displaystyle\sum_{j=1}^l\sum_{m=1}^j
  {\bf E}\left[{f_k}'''(S_{k-l-1})X_{k-m}\right]
     {\bf E}\left[X_{k-j}X_k\right]=$
\begin{eqnarray*}
&=&\sum_{j=1}^l\sum_{m=1}^j
  {\bf E}\left[{f_k}'''(S_{k-l-1})X_{k-m}\right]
     {\bf E}\left[X_{k-j}X_k\right]{\bf 1}_{\{l\ge(r+2)m\}}\\
&{ }& +\sum_{j=1}^l\sum_{m=1}^j
  {\bf E}\left[{f_k}'''(S_{k-l-1})X_{k-m}\right]
     {\bf E}\left[X_{k-j}X_k\right]{\bf 1}_{\{l< (r+2)m\}}\\
\end{eqnarray*}
Nous traitons la première des deux sommes du second membre :
$$\left\vert {\bf E}\left[{f_k}'''(S_{k-l-1})
     X_{k-m}\right]{\bf E}\left[X_{k-j}X_k\right]\right\vert
     \le CM^2{2(K_3+K_4)\over\sqrt{n}}
        {A\over n-k+\varepsilon^2}
        \varphi_{l+1-m,0}, $$
à l'aide de (\ref{sec:a2}) et de (\ref{sec:majo1}).
Il faut ensuite faire la somme sur les couples $(j,m)$ tels que 
$1\le m\le j \le l$
et $l\ge (r+2)m$.
Commençons  par remarquer que, si on a $l\ge (r+2)m$, 
alors on a $m\le
{\sqrt{n-k}\over r+2}$ et donc $l+1-m\ge {(r+1)\sqrt{n-k}\over r+2}$ et,
d'autre part, on a $l+1-m\ge (r+1)m$, donc $m\le {l+1-m\over r+1}$.
Ainsi, on a :
$$ \sum_{j=1}^l\sum_{m=1}^j\varphi_{l+1-m,m}{\bf 1}_{\{ l\ge (r+2)m\}}\le
\sum_{m=1}^{l}(l+1-m)\varphi_{l+1-m,m}{\bf 1}_{\{l\ge(r+2)m\}} $$
$$ \le \sum_{m=1}^{\left\lfloor{\sqrt{n-k}\over r+2}\right\rfloor}
  (l+1-m)\max_{j\le {l+1-m\over
r+1}}\varphi_{l+1-m,j}$$
$$\le \sum_{p=\left\lceil{(r+1)\sqrt{n-k}\over r+2}\right\rceil}^l p\max_{j\le{p\over
r+1}}\varphi_{p,j}\le\gamma_{\sqrt{n-k}},
$$
Puis traitons la seconde :
$$
\left\vert {\bf E}\left[{f_k}'''(S_{k-l-1})X_{k-m}
\right] {\bf E}\left[X_{k-j}X_k\right]\right\vert\le
      2{K_3\over\sqrt{n}}
      {A\over n-k+\varepsilon^2}M
         2CM\varphi_{j,0}
$$
Ainsi, on a :

$\displaystyle\sum_{j=1}^l\sum_{m=1}^j\left\vert
  {\bf E}\left[{f_k}'''(S_{k-m-1})X_{k-m}\right]
     {\bf E}\left[X_{k-j}X_k\right]\right\vert
     {\bf 1}_{\{l<(r+2)m\}}\le$
\begin{eqnarray}
&\leq& 
\sum_{j=1}^l\sum_{m=1}^j
 {4CM^2K_3\over\sqrt{n}}
      {A\over n-k+\varepsilon^2}
         \varphi_{j,0}{\bf 1}_{\{l<(r+2)m\}}\nonumber\\
&\leq&\sum_{j=\left\lceil{\sqrt{n-k}\over  {r+2}}\right\rceil}^l
  \sum_{m=1}^j{4CM^2K_3\over\sqrt{n}}
      {A\over n-k+\varepsilon^2}
         \varphi_{j,0}\nonumber\\
&\leq&{4CM^2K_3\over\sqrt{n}}
      {A\over n-k+\varepsilon^2}
      \sum_{j=\left\lceil{\sqrt{n-k}\over  {r+2}}\right\rceil}^l
  j\varphi_{j,0}\nonumber\\
&\leq&{4CM^2K_3\over\sqrt{n}}
      {A\over n-k+\varepsilon^2}
      \gamma_{\sqrt{n-k}}.
\end{eqnarray}

\item Reste le troisième terme

$\displaystyle\sum_{j=1}^l\sum_{m=1}^j\sum_{p=m+1}^{(r+2)j}\left\vert
   {\bf E}\left[\left({f_k}'''(S_{k-p})-{f_k}'''(S_{k-p-1})\right)
   X_{k-m}\right]{\bf E}\left[X_{k-j}X_k\right]\right\vert\le$
\begin{eqnarray}
&\le&\sum_{j=1}^l\sum_{m=1}^j\sum_{p=m+1}^{(r+2)j}
{2K_4M^2\over\sqrt{n}}{A\over(n-k+\varepsilon^2)^{3\over 2}}CM\varphi_{j,0}
\nonumber\\
&\le& 
   {2K_4M^2\over\sqrt{n}}{A\over (n-k+\varepsilon^2)^{3\over 2}}(r+2)
   CM\sum_{j=1}^{\lfloor \sqrt{n-k}\rfloor}j^2\varphi_{j,0}\nonumber\\
&\le& 
   {2CK_4(r+2)M^3\over\sqrt{n}}{A\over (n-k+\varepsilon^2)^{3\over 2}}(r+2)
   \alpha_{\sqrt{n-k}}
\end{eqnarray}
et

$\displaystyle\sum_{j=1}^l\sum_{m=1}^j\sum_{p=(r+2)j+1}^{l}\left\vert
   {\bf E}\left[\left({f_k}'''(S_{k-p})-{f_k}'''(S_{k-p-1})\right)
   X_{k-m}\right]{\bf E}\left[X_{k-j}X_k\right]\right\vert\le$
\begin{eqnarray}
&\le& \sum_{j=1}^l\sum_{m=1}^j\sum_{p=(r+2)j+1}^l
   {6K_4\over\sqrt{n}}{A\over (n-k+\varepsilon^2)^{3\over 2}}
   \varphi_{p-m,0}M^4\nonumber \\
&\le&{6K_4M^4\over\sqrt{n}}{A\over (n-k+\varepsilon^2)^{3\over 2}}
   \sum_{j=1}^l\sum_{m=1}^j\sum_{p=(r+2)j+1-m}^{l-m}\varphi_{p,0}\nonumber\\
&\le&{6K_4M^4\over\sqrt{n}}{A\over (n-k+\varepsilon^2)^{3\over 2}}
   \sum_{p=1}^l\sum_{j=1}^{\left\lfloor{p\over r+1}\right\rfloor}
    \sum_{m=1}^{j}\varphi_{p,0}\nonumber\\
&\le&{6K_4M^4\over\sqrt{n}}{A\over (n-k+\varepsilon^2)^{3\over 2}}
   {1\over(r+1)^2}\sum_{p=1}^l p^2\varphi_{p,0}\nonumber\\
&\le&{6K_4M^4\over\sqrt{n}}{1\over(r+1)^2}{A\alpha_{\sqrt{n-k}}
   \over (n-k+\varepsilon^2)^{3\over 2}} .
\end{eqnarray}

\end{enumerate}

\subsection{Contrôle de $\sum_{j=1}^lCov({f_k}''(S_{k-j-1}),X_{k-j}X_k)$}
Soit un entier $j=1,...,l$. Nous avons~:
$${f_k}''\left(S_{k-j-1}\right)={f_k}''\left(S_{k-l-1}\right)+
\sum_{m=j+1}^l\left({f_k}''\left(S_{k-m}\right)-
    {f_k}''\left(S_{k-m-1}\right)\right) $$
et, pour tout entier $m=j+1,...,l$~:
\begin{eqnarray*}
{f_k}''\left(S_{k-m}\right)-{f_k}''\left(S_{k-m-1}\right)& =&
     {f_k}'''(S_{k-m-1})X_{k-m}+ X_{k-m}^2\int_0^1(1-t)f^{(4)}\left(S_{k-m-1}
+tX_{k-m}\right)\, dt\\
&=& {f_k}'''(S_{k-l-1})X_{k-m}+\sum_{p=m+1}^l
\left({f_k}'''(S_{k-p})-{f_k}'''(S_{k-p-1})\right)X_{k-m}\\
&{ }& \ \ \ +X_{k-m}^2\int_0^1(1-t)f^{(4)}\left(S_{k-m-1}
+tX_{k-m}\right)\, dt.
\end{eqnarray*}
Nous pouvons donc écrire :
\begin{eqnarray*}
\sum_{j=1}^lCov({f_k}''(S_{k-j-1}),X_{k-j}X_k)&=&\sum_{j=1}^lCov( {f_k}''\left(S_{k-l-1}\right),X_{k-j}X_k)\\
&+& \sum_{j=1}^l\sum_{m=j+1}^l
Cov\left(X_{k-m}^2\int_0^1(1-t)f^{(4)}\left(S_{k-m-1}
+tX_{k-m}\right)\, dt,X_{k-j}X_k\right)\\
&+& \sum_{j=1}^l\sum_{m=j+1}^lCov( {f_k}'''(S_{k-l-1})X_{k-m},X_{k-j}X_k)\\
&+& \sum_{j=1}^l\sum_{m=j+1}^lCov\left(\sum_{p=m+1}^l
\left({f_k}'''(S_{k-p})-{f_k}'''(S_{k-p-1})\right)X_{k-m},X_{k-j}X_k\right),\\
\end{eqnarray*}
et, finalement~:
\begin{eqnarray*}
\sum_{j=1}^lCov({f_k}''(S_{k-j-1}),X_{k-j}X_k)&=&\sum_{j=1}^l
Cov\left( {f_k}''\left(S_{k-l-1}\right),X_{k-j}X_k\right)\\
&+& \sum_{j=1}^lCov(X_{k-m}^2\int_0^1(1-t)f^{(4)}\left(S_{k-m-1}
+tX_{k-m}\right)\, dt,X_{k-j}X_k)\\
&+& \sum_{j=1}^l\sum_{m=j+1}^lCov( {f_k}'''(S_{k-l-1})
X_{k-m},X_{k-j}X_k){\bf 1}_{\{l<(r+2)m\}}\\
&+& \sum_{j=1}^l\sum_{m=j+1}^lCov( {f_k}'''(S_{k-l-1})
X_{k-m},X_{k-j}X_k){\bf 1}_{\{l\ge(r+2)m\}}\\
&+& \sum_{j=1}^l\sum_{m=j+1}^l\sum_{p=m+1}^{\min(l,(r+2)m)}
Cov\left(\left({f_k}'''(S_{k-p})-{f_k}'''(S_{k-p-1})\right)X_{k-m},X_{k-j}
X_k\right)\\
&+& \sum_{j=1}^l\sum_{m=j+1}^l
Cov\left(\left({f_k}'''(S_{k-l-1})-{f_k}'''(S_{k-(r+2)m-1})
\right)X_{k-m},X_{k-j}
X_k\right){\bf 1}_{\{l\ge(r+2)m+1\}}.
\end{eqnarray*}
Nous allons majorer chacun des termes apparaissant dans le second membre 
de cette égalité.
\begin{enumerate}
\item {\bf (Contrôle de 
$\sum_{j=1}^lCov( {f_k}''\left(S_{k-l-1}\right),X_{k-j}X_k)$)}

\'Ecrivons :
\begin{eqnarray*}
\sum_{j=1}^lCov( {f_k}''\left(S_{k-l-1}\right),X_{k-j}X_k)&=&
\sum_{j=1}^{\left\lfloor{l\over r+1}\right\rfloor}Cov( {f_k}''\left(S_{k-l-1}\right),X_{k-j}X_k)\\
&{ }&\ \ \ \ \ +\sum_{j=\left\lfloor{l\over r+1}\right\rfloor+1}
   ^lCov( {f_k}''\left(S_{k-l-1}\right),X_{k-j}X_k).\\
\end{eqnarray*}
D'après (\ref{sec:a2}) et (\ref{sec:majo1}), nous avons~:

$\displaystyle
\left\vert Cov( {f_k}''\left(S_{k-l-1}\right),X_{k-j}X_k)\right\vert$
\begin{eqnarray*}
&\le& \left\vert Cov( {f_k}''\left(S_{k-l-1}\right)X_{k-j},X_k)\right\vert
   +\left\vert {\bf E}\left[ {f_k}''\left(S_{k-l-1}\right)\right]
          {\bf E}\left[X_{k-j}X_k\right]\right\vert\\
&\le& {8C(K_2+K_3)M\over\sqrt{n}}
      {A\over\sqrt{ n-k+\varepsilon^2}}\varphi_{j,0}
\end{eqnarray*}
et 
$$ 
\left\vert
Cov( {f_k}''\left(S_{k-l-1}\right),X_{k-j}X_k)\right\vert\le 
  {2C(K_2+K_3)\over\sqrt{n}}
      {A\over\sqrt{ n-k+\varepsilon^2}}\varphi_{l-j,j}.
$$
En reportant ces majorations dans chacune des deux sommes on obtient :

$\displaystyle \sum_{j=1}^l
\left\vert Cov( {f_k}''\left(S_{k-l-1}\right),X_{k-j}X_k)\right\vert\le$
\begin{eqnarray}
&\le& 
{8C(K_2+K_3)M\over\sqrt{n}}
      {A\over\sqrt{ n-k+\varepsilon^2}}
      \left(\sum_{j=1}^{\left\lfloor{l\over r+1}\right\rfloor}
      \varphi_{l-j,j}+
      \sum_{j=\left\lfloor{l\over r+1}\right\rfloor+1}^l\varphi_{j,0}\right)
      \nonumber\\
      &\le& 
{8C(K_2+K_3)M\over\sqrt{n}}
      {A\over\sqrt{ n-k+\varepsilon^2}}
      \delta_{n,k}.
\end{eqnarray}

\item {\bf (Contrôle de $\sum_{j=1}^l\sum_{m=j+1}^l
    Cov\left(X_{k-m}^2\int_0^1(1-t)f^{(4)}\left(S_{k-m-1}
+tX_{k-m}\right)\, dt,X_{k-j}X_k\right)$)}
   
Soit un couple d'entiers $(j,m)$ vérifiant~: $1\le j\le j+1\le m\le l$.
\begin{itemize}
\item Lorsque $m\le (r+2)j$, on a~:

$\displaystyle\left\vert Cov\left(X_{k-m}^2\int_0^1(1-t)f^{(4)}\left(S_{k-m-1}
+tX_{k-m}\right)\, dt\ X_{k-j},X_k\right)\right\vert=$
\begin{eqnarray*}
&=&\left\vert \int_0^1(1-t)
 Cov\left(f^{(4)}\left(S_{k-m-1}
+tX_{k-m}\right)X_{k-m}^2 X_{k-j},X_k\right)\, dt\right\vert\\
&\le& C{18(K_4+K_5)M^3\over\sqrt n}{A\over(n-k+\varepsilon^2)^{3\over 2}}
 \varphi_{j,0}.
\end{eqnarray*}
 et
$$\left\vert {\bf E}\left[X_{k-m}^2\int_0^1(1-t)f^{(4)}\left(S_{k-m-1}
+tX_{k-m}\right)\, dt \right]{\bf E}\left[X_{k-j}X_k
    \right]\right\vert
 \le C{K_4\over \sqrt{n}}{A\over (n-k+\varepsilon^2)^{3\over 2}}M^3
 \varphi_{j,0} . $$
\item Lorsque $m\ge (r+2)j+1$, alors on a~:
\begin{eqnarray*}
&{ }&  \left\vert Cov\left(X_{k-m}^2\int_0^1(1-t)f^{(4)}\left(S_{k-m-1}
+tX_{k-m}\right)\, dt ,X_{k-j}X_k\right)\right\vert\\
&{ }& \ \ \ \ \ \
 \le C{10(K_4+K_5)M^2\over\sqrt{n}}
 {A\over(n-k+\varepsilon^2)^{3\over 2}}\varphi_{m-j,j}. 
\end{eqnarray*}
\end{itemize}
Ainsi, en utilisant la relation 
$$
Cov(A,BC)=Cov(AB,C)-{\bf E}(A){\bf E}(BC)
$$
valable quand ${\bf E}(C)$ est nulle, 
on a~:

$\displaystyle\sum_{j=1}^l\sum_{m=j+1}^l
    \left\vert Cov\left(X_{k-m}^2\int_0^1(1-t)f^{(4)}\left(S_{k-m-1}
+tX_{k-m}\right)\, dt ,X_{k-j}X_k\right)\right\vert\le$
\begin{eqnarray}
&\le& {\tilde K_0\over\sqrt{n}}{A\over(n-k+\varepsilon^2)^{3\over 2}}
      \left(\sum_{p=1}^{\lfloor \sqrt{n-k}\rfloor}p\varphi_{p,0}
      +\sum_{p=1}^{\lfloor\sqrt{n-k}
        \rfloor}\sum_{l=1
      }^{\left\lfloor{p\over r+1}\right\rfloor}\varphi_{p,l}\right)
      \nonumber\\
&\le& {\tilde K_0\over\sqrt{n}}{A\alpha_{\sqrt{n-k}}
\over(n-k+\varepsilon^2)^{3\over 2}}.
\end{eqnarray}
pour un certain $\tilde K_0$ dépendant de $C$, $K_4$, $K_5$, $M$ et $r$.
\item Nous avons :

$\displaystyle\sum_{j=1}^l\sum_{m=j+1}^{l}
\left\vert Cov\left({f_k}'''(S_{k-l-1})
    X_{k-m},X_{k-j}X_k\right){\bf 1}_{\{l< (r+2)m\}}\right\vert\le$
\begin{eqnarray}
&\leq & 
    {\tilde K_1\over\sqrt{n}}{A\over n-k+\varepsilon^2}
     \left( \sum_{p=\left\lceil{\sqrt{n-k}\over (r+2)^2}\right\rceil}^
      {\lfloor\sqrt{n-k}\rfloor}p\varphi_{p,0}
     +\sum_{p=\left\lceil{(r+1)\sqrt{n-k}\over (r+2)^2}\right\rceil}^{
       \lfloor \sqrt{n-k}\rfloor}\sum_{j=1}
       ^{\left\lfloor {p\over r+1}\right\rfloor}\varphi_{p,j}\right)\nonumber\\
&\leq & 
    {\tilde K_1\over\sqrt{n}}{A\gamma_{\sqrt{n-k}}\over n-k+\varepsilon^2},
\end{eqnarray}
pour un certain $\tilde K_1$ dépendant de $C$, $K_3$, $K_4$, $M$ et $r$.
En effet, considérons un couple d'entiers $(j,m)$ vérifiant~:
$1\le j\le j+1\le m\le l$.
\begin{itemize}
\item Si on a $l< (r+2)m$ et $m\le (r+2)j$, alors on a~:
$$\left\vert Cov\left({f_k}'''(S_{k-l-1})
X_{k-m}X_{k-j},X_k\right)\right\vert\le
      C{8(K_3+K_4)\over\sqrt{n}}
   {A\over n-k+\varepsilon^2}
       M^2\varphi_{j,0} $$
et
$$\left\vert {\bf E}\left[{f_k}'''(S_{k-l-1})X_{k-m}
\right] {\bf E}\left[X_{k-j}X_k\right]\right\vert\le
      C{2K_3M\over\sqrt{n}}
      {A\over n-k+\varepsilon^2}
         2M\varphi_{j,0} .$$
\item Si on a $l< (r+2)m$ et $m> (r+2)j$, alors on a~:
$$\left\vert Cov\left({f_k}'''(S_{k-l-1})
X_{k-m},X_{k-j}X_k\right)\right\vert\le
      C{4(K_3+K_4)M\over\sqrt{n}}
      {A\over n-k+\varepsilon^2}
         \varphi_{m-j,j} .$$
\end{itemize}

\item Nous avons~:

$\displaystyle\sum_{j=1}^l\sum_{m=j+1}^{l} Cov\left({f_k}'''
(S_{k-l-1})
    X_{k-m},X_{k-j}X_k\right){\bf 1}_{\{l\ge(r+2)m\}}=$
\begin{eqnarray*}
&=& \sum_{j=1}^l\sum_{m=j+1}^{l}{\bf E}[{f_k}'''(S_{k-l-1})]
{\bf E}[ X_{k-m}X_{k-j}X_k]{\bf 1}_{\{l\ge(r+2)m\}}\\
&+& \sum_{j=1}^l\sum_{m=j+1}^{l} Cov({f_k}'''(S_{k-l-1}), X_{k-m}X_{k-j}X_k)
{\bf 1}_{\{l\ge(r+2)m\}}\\
&+& \sum_{j=1}^l\sum_{m=j+1}^{l} {\bf E}[{f_k}'''(S_{k-l-1})X_{k-m}]
{\bf E}[X_{k-j}X_k]{\bf 1}_{\{l\ge(r+2)m\}}.\\
\end{eqnarray*}
Nous allons contrôler chaque terme du membre de droite de cette identité.
\begin{itemize}
\item Nous commençons par majorer le premier terme~:

$\displaystyle\sum_{j=1}^l\sum_{m=j+1}^{l}\left\vert{\bf E}
[{f_k}'''(S_{k-l-1})]
{\bf E}[ X_{k-m}X_{k-j}X_k]\right\vert\le$
$$\le K_3\left({A\over(n-k+\varepsilon^2)^{3\over 2}}+{1\over n}\right)
\sum_{j=1}^l\sum_{m=j+1}^{l}\left\vert{\bf E}[ X_{k-m}X_{k-j}X_k]\right\vert.
$$
Or, nous avons~:
\begin{eqnarray*}
\sum_{j=1}^l\sum_{m=j+1}^{(r+2)j}\left\vert{\bf E}[ X_{k-m}X_{k-j}X_k]
\right\vert
&=&\sum_{j=1}^l\sum_{m=j+1}^{(r+2)j}\left\vert Cov( X_{k-m}X_{k-j},X_k)
  \right\vert\\
&\le&\sum_{j=1}^l\sum_{m=j+1}^{(r+2)j}3CM^2\varphi_{j,0}\\
 &\le&3CM^2\sum_{j=1}^l (r+1)j\varphi_{j,0}
\end{eqnarray*}
et
\begin{eqnarray*}
\sum_{j=1}^l\sum_{m=(r+2)j+1}^{l}
\left\vert{\bf E}[ X_{k-m}X_{k-j}X_k]\right\vert
&=&\sum_{j=1}^l\sum_{m=(r+2)j+1}^{l}\left\vert Cov( X_{k-m,}X_{k-j}X_k)
  \right\vert\\
&\le&\sum_{j=1}^l\sum_{m=(r+2)j+1}^{l}  2CM\varphi_{m-j,j}\\
   &\le& 
2CM\sum_{p=1}^l\sum_{j=1}^{\left\lfloor{p\over r+1}\right\rfloor} 
 \varphi_{p,j}.
\end{eqnarray*}
D'où~:
\begin{equation}
\sum_{j=1}^l\sum_{m=j+1}^{l}\left\vert{\bf E}
[{f_k}'''(S_{k-l-1})]
{\bf E}[ X_{k-m}X_{k-j}X_k]\right\vert\le
5(r+1)CM^2K_3\left({A\alpha_{\sqrt{n-k}}\over(n-k+\varepsilon^2)^{3\over 2}}+{
   \beta_{\sqrt{n-k}}\over n}\right).
\end{equation}
\item Pour les deux autres termes,  il suffit d'écrire :
$$
\left\vert  Cov\left({f_k}'''(S_{k-l-1}),
     X_{k-m}X_{k-j}X_k\right)
     \right\vert\le C{2(K_3+K_4)\over\sqrt{n}}
        {A\over n-k+\varepsilon^2}\varphi_{l+1-m,m}$$
et
$$\left\vert {\bf E}\left[{f_k}'''(S_{k-l-1})
     X_{k-m}\right]{\bf E}\left[X_{k-j}X_k\right]\right\vert
     \le CM^2{2(K_3+K_4)\over\sqrt{n}}
        {A\over n-k+\varepsilon^2}
        \varphi_{l+1-m,0}, $$
à l'aide de la propriété de décorrélation (\ref{sec:a2}) et de 
l'inégalité (\ref{sec:majo1}).
Il faut ensuite faire la somme sur les couples $(j,m)$ avec $1\le j<m\le l$
tels que $l\ge (r+2)m$.
Commençons  par remarquer que, si on a $l\ge (r+2)m$; alors on a 
$m\le
{\sqrt{n-k}\over r+2}$ et donc $l+1-m\ge {(r+1)\sqrt{n-k}\over r+2}$;
d'autre part on a $l+1-m\ge (r+1)m$, donc $m\le {l+1-m\over r+1}$.
Ainsi, on a :
$$ \sum_{j=1}^l\sum_{m=j+1}^l\varphi_{l+1-m,m}{\bf 1}_{\{ l\ge (r+2)m\}}\le
\sum_{m=2}^{l}(m-1)\varphi_{l+1-m,m}{\bf 1}_{\{l\ge(r+2)m\}} $$
$$ \le \sum_{m=2}^{\left\lfloor{\sqrt{n-k}\over r+1}\right\rfloor}
   {l+1-m\over r+1}\max_{j\le {l+1-m\over
r+1}}\varphi_{l+1-m,j}$$
$$\le \sum_{p=\left\lceil{(r+1)\sqrt{n-k}\over r+2}\right\rceil}^l 
p\max_{j\le{p\over
r+1}}\varphi_{p,j}\le\gamma_l.$$
Nous obtenons ainsi~:
\begin{equation}
\sum_{j=1}^l\sum_{m=j+1}^{l} \left\vert
Cov({f_k}'''(S_{k-l-1}), X_{k-m}X_{k-j}X_k)\right\vert
{\bf 1}_{\{l\ge(r+2)m\}}\le C{2(K_3+K_4)\over\sqrt{n}}
        {A\gamma_{\sqrt{n-k}}\over n-k+\varepsilon^2}
\end{equation}
De même, nous obtenons~:
\begin{equation}
\sum_{j=1}^l\sum_{m=j+1}^{l} \left\vert{\bf E}[{f_k}'''(S_{k-l-1})X_{k-m}]
{\bf E}[X_{k-j}X_k]\right\vert{\bf 1}_{\{l\ge(r+2)m\}}
 \le CM^2{2(K_3+K_4)\over\sqrt{n}}
        {A\gamma_{\sqrt{n-k}}\over n-k+\varepsilon^2}.
\end{equation}
\end{itemize}
\item Nous avons~:

$\displaystyle\sum_{j=1}^l\sum_{m=j+1}^{l}\sum_{p=m+1}^{\min(l,(r+2)m)}
\left\vert Cov\left(\left({f_k}'''(S_{k-p})-{f_k}'''(S_{k-p-1})
\right)X_{k-m},X_{k-j}X_k\right)\right\vert\le $
\begin{eqnarray}
&\le& 
    {\tilde K_3\over\sqrt{n}}{A\over (n-k+\varepsilon^2)^{3\over 2}}
     \left( \sum_{p=1}^{\lfloor\sqrt{n-k}\rfloor}p^2\varphi_{p,0}
     +\sum_{p=1}^{
       \lfloor \sqrt{n-k}\rfloor}p\sum_{j=1}
       ^{\left\lfloor {p\over r+1}\right\rfloor}\varphi_{p,j}\right)\nonumber\\
&\le& {\tilde K'_3\over\sqrt{n}}{A\alpha_{\sqrt{n-k}}
\over (n-k+\varepsilon^2)^{3\over 2}},
   \label{sec:A}
\end{eqnarray}
pour un certain $\tilde K_3$ dépendant de $C$, $K_4$, $M$ et $r$.
En effet, considérons un triplet d'entiers $(j,m,p)$ vérifiant~:
$1\le j\le j+1\le m\le m+1\le p\le l$.
\begin{itemize}
\item Si on a $p\le (r+2)m$ et $m\le (r+2)j$, alors on a~:
\begin{equation}
\left\vert Cov\left(\left({f_k}'''(S_{k-p})-{f_k}'''(S_{k-p-1}\right)
X_{k-m}X_{k-j},X_k\right)\right\vert\le
      C{18K_4\over\sqrt{n}}
       {A\over (n-k+\varepsilon^2)^{3\over 2}}M^3
          \varphi_{j,0} 
\end{equation}
et
\begin{equation}
\left\vert {\bf E}\left[\left({f_k}'''(S_{k-p})-
{f_k}'''(S_{k-p-1})\right)X_{k-m}
\right] {\bf E}\left[X_{k-j}X_k\right]\right\vert\le
      C{4K_4\over\sqrt{n}}
       {A\over (n-k+\varepsilon^2)^{3\over 2}}M^3\varphi_{j,0} .
\end{equation}
\item D'autre part, si on a $p\le (r+2)m$ et $m\ge (r+2)j+1$, alors on a~:
$$
\left\vert Cov\left(\left({f_k}'''(S_{k-p})-{f_k}'''(S_{k-p-1})\right)
X_{k-m},X_{k-j}X_k\right)\right\vert\le
      C{6K_4\over\sqrt{n}}
       {A\over (n-k+\varepsilon^2)^{3\over 2}}M^2M\varphi_{m-j,j} .$$
\end{itemize}
Nous concluons en remarquant que~:
$$\displaystyle \sum_{j=1}^l\sum_{m=j+1}^{(r+2)j}\sum_{p=m+1}^{(r+2)m}
   \varphi_{j,0}\le(r+2)^2\sum_{j=1}^lj^2\varphi_{j,0}$$
et que
\begin{eqnarray*}
\sum_{j=1}^l\sum_{m=(r+2)j+1}^{l}\sum_{p=m+1}^{(r+2)m}\varphi_{m-j,j}
&=&\sum_{j=1}^l\sum_{m=(r+2)j+1}^{l}(r+1)m\varphi_{m-j,j}\\
&=&\sum_{j=1}^l\sum_{p=(r+1)j+1}^{l-j}(r+1)(p+j)\varphi_{p,j}\\
&\le&\sum_{p=1}^l\sum_{j=1}^{\left\lfloor{p\over r+1}\right\rfloor}(r+1)(p+j)
\varphi_{p,j}\\
&\le&(r+2)
\sum_{p=1}^l\sum_{j=1}^{\left\lfloor{p\over r+1}\right\rfloor}p
\varphi_{p,j}.
\end{eqnarray*}

\item Nous contrôlons à présent~:
$$ \sum_{j=1}^l\sum_{m=j+1}^l
Cov\left(\left({f_k}'''(S_{k-l-1})-{f_k}'''(S_{k-(r+2)m-1})
\right)X_{k-m},X_{k-j}
X_k\right){\bf 1}_{\{l\ge(r+2)m+1\}}.$$
Si $l\ge (r+2)m+1$, nous avons~:
$$
{f_k}'''(S_{k-l-1})-{f_k}'''(S_{k-(r+2)m-1})=
   \sum_{p=(r+2)m+1}^{l}{f_k}'''(S_{k-p})-{f_k}'''(S_{k-p-1}).
$$
Nous allons utiliser l'identité suivante~:
$$
Cov(AB,CD)=Cov(A,BCD)-{\bf E}[AB]{\bf E}[CD]+{\bf E}[A]{\bf E}[BCD].
$$

\begin{itemize}
\item 
%
Quand $l\ge p\ge (r+2)m+1$, on a~:
$$\left\vert  Cov\left( {f_k}'''(S_{k-p})-{f_k}'''(S_{k-p-1}),
    X_{k-m}X_{k-j}X_k\right)
    \right\vert\le C{6K_4
    A M\over{\sqrt{n}}(n-k+\varepsilon^2)^{3\over 2}}
    \varphi_{p-m,m}$$
et
$$\left\vert {\bf E}\left[{f_k}'''(S_{k-p})-{f_k}'''(S_{k-p-1})
    X_{k-m}\right]{\bf E}\left[X_{k-j}X_k\right]\right\vert
    \le C{6K_4AM^3\over{\sqrt{n}}(n-k+\varepsilon^2)^{3\over 2}}
    \varphi_{p-m,0}. 
$$
Nous majorons donc les sommes de ces quantités sur les couples d'entiers
$(j,m,p)$ vérifiant $1\le j\le j+1\le m$ et $ (r+2)m+1\le p\le l$ par~:
$$
{\tilde K_4\over\sqrt{n}}{A\over (n-k+\varepsilon^2)^{3\over 2}}
     \left( \sum_{p=1}^{\lfloor\sqrt{n-k}\rfloor}p^2\varphi_{p,0}
     +\sum_{p=1}^{
       \lfloor \sqrt{n-k}\rfloor}p\sum_{j=1}
       ^{\left\lfloor {p\over r+1}\right\rfloor}\varphi_{p,j}\right).
$$
Nous avons donc~:

$\displaystyle\sum_{j=1}^l\sum_{m=j+1}^l\left\vert
Cov\left({f_k}'''(S_{k-l-1})-{f_k}'''(S_{k-(r+2)m-1})
,X_{k-m}X_{k-j}
X_k\right)\right\vert{\bf 1}_{\{l\ge(r+2)m+1\}}\le $
\begin{equation}
\le{\tilde K_4\over\sqrt{n}}{A\alpha_{\sqrt{n-k}}
\over (n-k+\varepsilon^2)^{3\over 2}}
\end{equation}
et

$\displaystyle\sum_{j=1}^l\sum_{m=j+1}^l\left\vert
{\bf E}\left[\left({f_k}'''(S_{k-l-1})-{f_k}'''(S_{k-(r+2)m-1})
\right)X_{k-m}\right]{\bf E}\left[X_{k-j}
X_k\right]\right\vert{\bf 1}_{\{l\ge(r+2)m+1\}}\le $
\begin{equation}
\le{\tilde K_4\over\sqrt{n}}{A\alpha_{\sqrt{n-k}}
\over (n-k+\varepsilon^2)^{3\over 2}}
\end{equation}
\item Contrôlons à présent la quantité suivante~:
$$
\sum_{j=1}^l\sum_{m=j+1}^{l}{\bf E}\left[{f_k}'''
(S_{k-l-1})-{f_k}'''(S_{k-(r+2)m-1})\right]{\bf E}
\left[X_{k-m}X_{k-j}X_k\right].
$$
Lorsque $l\ge (r+2)m+1$ et $m\ge(r+2)j+1$,  on a~:
$$
\left\vert {\bf E}\left[{f_k}'''(S_{k-(r+2)m-1})
\right] Cov\left(X_{k-m},X_{k-j}X_k\right)\right\vert\le
      {K_3\over\sqrt{n}}
      \left({A\over (n-k+\varepsilon^2)^{3\over 2}}
      +{1\over n}\right)2CM\varphi_{m-j,j} .$$
Quand $l\ge  (r+2)m+1$ et si $m\le (r+2)j$, on a~:
$$\left\vert {\bf E}\left[{f_k}'''(S_{k-(r+2)m-1})
\right] Cov\left(X_{k-m}X_{k-j},X_k\right)\right\vert\le
{K_3\over\sqrt{n}}
      \left({A\over (n-k+\varepsilon^2)^{3\over 2}}
      +{1\over n}\right)C3M^2\varphi_{j,0}.$$
On en déduit :

$\displaystyle 
\sum_{j=1}^l\sum_{m=j+1}^l\left\vert {\bf E}\left[f_k'''
(S_{k-(r+2)m-1})
\right] {\bf E}\left[X_{k-m}X_{k-j}X_k\right]{\bf 1}_{\{l\ge(r+2)m+1\}}
\right\vert\le$
\begin{eqnarray}
  &\le& {\tilde K_6\over\sqrt{n}}
  \left({A\over (n-k+\varepsilon^2)^{3\over 2}}
  +{1\over n}
  \right)\left(\sum_{p=1}^{
     \lfloor\sqrt{n-k}\rfloor}\sum_{j=1}^{\lfloor {p\over r+1}
      \rfloor}\varphi_{p,j}+\sum_{p=1}^{\lfloor \sqrt{n-k}
      \rfloor}p\varphi_{p,0}\right) \nonumber\\
&\le&{\tilde K_6\over\sqrt{n}}
  \left({A\alpha_{\sqrt{n-k}}\over (n-k+\varepsilon^2)^{3\over 2}}
  +{\beta_{\sqrt{n-k}}\over n}
  \right)
\end{eqnarray}
De la même façon, on montre la majoration suivante :

$
\displaystyle 
\sum_{j=1}^l\sum_{m=j+1}^l\left\vert {\bf E}\left[f_k'''(S_{k-l-1})
\right] {\bf E}\left[X_{k-m}X_{k-j}X_k\right]{\bf 1}_{\{l\ge(r+2)m+1\}}
\right\vert\le
$
\begin{equation}
  \le{\tilde K_6\over\sqrt{n}}
  \left({A\alpha_{\sqrt{n-k}}\over (n-k+\varepsilon^2)^{3\over 2}}
  +{\beta_{\sqrt{n-k}}\over n}
  \right).
\end{equation}

\end{itemize}
\end{enumerate}

\end{document}